\documentclass[10pt,twoside]{article}
\usepackage{lmodern}
\usepackage{cite}
\usepackage[T1]{fontenc}
\usepackage{amsthm}
\usepackage{amsfonts}
\usepackage{amssymb}
\usepackage{amsmath}
\usepackage{graphicx}
\usepackage{mathrsfs} % for \mathscr{F}
\usepackage{mathtools}
\usepackage{hyperref}
\usepackage{array,multirow}
\usepackage{caption}   % for tables
\usepackage{float}   % for tables
\usepackage{microtype}
\usepackage{subfiles}    %for splitting input file

\usepackage{tikz}
\usetikzlibrary{matrix}
\usetikzlibrary{calc} 

   \author{by  Martin~Seysen\\
}
   \title{A fast implementation of the Monster group \\[0.3em]
   	\Large	
   	The Monster has been tamed
   }

   \newcommand{\Field}{\mathbb{F}}
   \newcommand{\setZ}{\mathbb{Z}}
   \newcommand{\setR}{\mathbb{R}}

   \newcommand{\setF}{\mathbb{F}}

   \newcommand{\setM}{\mathbb{M}}

   \newcommand{\Cfrac}[1]{{/ \kern-0.30em / #1 / \kern-0.30em /}}

   \newcommand{\Skip}[1]{}

   \newcommand {\End}{\,{\mathop{\rm E\kern-0.05em{n}\kern-0.05em{d}}}}

   \newcommand {\rad}{\,{\mathop{\rm  rad}}}

  \newcommand {\Aut}{\,{\mathop{\rm Aut}}\,}

 \newcommand {\Span}{\,{\mathop{\rm span}}}
 \newcommand {\rank}{\,{\mathop{\rm rank}}}

   \newcommand{\AutStP}{\Aut_{\!\mbox{\scriptsize St}} \mathcal{P}}

\DeclareMathSymbol{\mlq}{\mathord}{operators}{``}
\DeclareMathSymbol{\mrq}{\mathord}{operators}{`'}

   \newcommand {\fatline}[1]{\noindent {\bf #1}}

   \newcommand {\proofend}{\noindent $\Box$ }

   \newtheorem {Theorem}{Theorem}[section]

   \newtheorem {Lemma}[Theorem]{Lemma}

   \numberwithin{equation}{Theorem}
 
   \newcounter{StatementListCounter}

\begin{document} {\large}

  \maketitle

  \begin{abstract}
  	
Let $\mathbb{M}$ be the  Monster group, which is the largest sporadic
finite simple  group, and has first been constructed in 1982 by 
Griess. In 1985 Conway has constructed a 196884-dimensional
rational representation $\rho$ of  $\mathbb{M}$ with matrix entries
in $\mathbb{Z}[\frac{1}{2}]$. We describe a new and very  fast 
algorithm for performing the group operation  in  $\mathbb{M}$.

For an odd integer $p > 1$ let $\rho_p$ be the representation
$\rho$ with matrix entries taken modulo $p$. 
We use a generating set $\Gamma$ of $\mathbb{M}$, such that the operation
of a generator in  $\Gamma$ on an element of $\rho_p$ can easily
be computed. 

We construct a triple  $(v_1, v^+, v^-)$ of elements of the module 
$\rho_{15}$, such that an unknown $g \in \mathbb{M}$ can be effectively
computed as a word in $\Gamma$ from the images  
$(v_1  g, v^+  g, v^- g)$. 

Our new algorithm based on this idea multiplies two random elements 
of  $\mathbb{M}$ in less than  30~milliseconds on a standard PC with
an Intel  i7-8750H CPU at 4 GHz. This is more than 100000 times faster 
than estimated by Wilson in 2013.

  \end{abstract}

\section*{}

{\noindent \bf Key Words:}

Monster group, finite simple groups, group representation, 
efficient implementation

\vspace{1ex}

\noindent MSC2020-Mathematics Subject Classification (2020):
20C34, 20D08, 20C11, 20--08

%%%%%%%%%%%%%%%%%%%%%%%%%%%%%%%%%%%%%%%%%%%%%%%%%%%%%%%%%%%%%%%%%%%%%%%%%%%%%%
\section{Introduction}
\label{sect:Introduction}
%%%%%%%%%%%%%%%%%%%%%%%%%%%%%%%%%%%%%%%%%%%%%%%%%%%%%%%%%%%%%%%%%%%%%%%%%%%%%%

Let $\setM$ be the  Monster group, which is the largest 
sporadic finite simple  group, and  has first been constructed
by Griess \cite{Griess:Friendly:Giant}. That construction has
been simplified by Conway \cite{Conway:Construct:Monster},
leading to a rational representation $\rho$ of $\setM$ of 
dimension 196884 with matrix entries in  $\setZ[\frac{1}{2}]$.
For a small odd integer $p > 1$ let  $\rho_p$ be the representation
$\rho$ with matrix entries taken modulo $p$. In this paper
we deal with $\rho_{15}$ and, occasionally, with  $\rho_{3}$
and  $\rho_{5}$.  

The first computer construction of the Monster group is due to
Linton, Parker, Walsh and Wilson \cite{LPWW-Computer_Monster}. That 
construction is based on large 3-local subgroups. The reason
for choosing 3-local subgroups (instead of the much larger 2-local 
subgroups) was that here computations can be done with scalars in 
$\mathbb{F}_2$.  Holmes and Wilson \cite{holmes_wilson_2003}
have presented a computer construction of the Monster based
on 2-local subgroups with scalars in $\mathbb{F}_3$. The 
representation used in  that construction resembles the 
representation $\rho_3$ mentioned above.
Here the 3-local construction appears to be faster,
but the 2-local construction allows general computations in
a much larger maximal 2-local subgroup of $\mathbb{M}$.

Any known faithful representation of $\mathbb{M}$ has dimension at 
least 196882; so in practice we cannot store elements of $\mathbb{M}$
as matrices acting on such a representation on a standard PC. We store 
elements of $\mathbb{M}$ as words in a set $\Gamma$ of  generators
of $\mathbb{M}$, where each generator corresponds to a sparse matrix
acting on $\rho$. We write $\Gamma^*$ for the set of words in
$\Gamma$. Inverting a word in $\Gamma^*$ is easy by construction
of $\Gamma$. Group multiplication in $\setM$ is simply a 
concatenation of words in $\Gamma^*$. But here the word shortening 
problem arises, since the length of a word may grow exponentially
with the number of group operations. Wilson \cite{Wilson13} has
presented a general word shortening algorithm for the Monster 
group.

In this paper we present a construction of the Monster based
on the representation $\rho$ mentioned above. We also give a new
word shortening algorithm. We construct a triple  $(v_1, v^+, v^-)$ 
of elements of $\rho_{15}$, such that for every
$g \in \setM$ we can effectively compute a word $g' \in \Gamma^*$ 
with $g' = g$ from the images  $(v_1  g, v^+  g, v^- g)$, without
using $g$. Obviously, the word $g'$ computed from these three images 
depends on the element $g$ of the Monster only, but not on the 
representation of $g$ as a word in $\Gamma^*$.

This means that we may effectively compute a unique {\em reduced} 
form of each element of the Monster as a word in
$\Gamma^*$. Using standard data compression methods, we may store
any such reduced form in less than 256 bits. So we may quickly 
find an element of the Monster in a list of millions of such 
elements. This was not possible before. 

We have implemented the group operation
in $\mathbb{M}$ using that new reduction algorithm in the software 
project~\cite{mmgroup2020}. For documentation of the project, see
\cite{mmgroup_doc}. The run time for the group operation in that
project is a bit less than 30 milliseconds on the author's PC,
which has an Intel i7-8750H CPU at about 4 GHz. 
So this is the first implementation of the group operation
of the Monster that a user can run interactively on a computer.

We have also implemented the word shortening algorithm 
in \cite{Wilson13}; here a group operation using that algorithm 
takes about 30 seconds. In 2013 the run time of this operation 
has been estimated to take 1--2 hours, see \cite{Wilson13}.

Dietrich, Lee, and Popiel \cite{DLP2023} have used our new
algorithm for solving the long-standing problem of finding all 
maximal subgroups of the Monster.

%%%%%%%%%%%%%%%%%%%%%%%%%%%%%%%%%%%%%%%%%%%%%%%%%%%%%%%%%%%%%%%%%%%%%%%%%%%%%%
\section{Overview of the new algorithm}
\label{sect:Overview}
%%%%%%%%%%%%%%%%%%%%%%%%%%%%%%%%%%%%%%%%%%%%%%%%%%%%%%%%%%%%%%%%%%%%%%%%%%%%%%

In this section we give a brief overview of the new reduction 
algorithm. More details will be given in the following sections.

\subsection{Construction of the Monster and of its representation}
\label{ssect:Construct:Monster}

In this subssection we briefly recap Conway's construction 
\cite{Conway:Construct:Monster} of the Monster.
A more detailed description of this construction given in 
Sections~\ref{sect:N0}~--~\ref{sect:rep:M}.

The Monster group contains two classes of involutions called 2A and 
2B in the ATLAS \cite{Atlas}. The construction of the Monster
$\setM$  in \cite{Conway:Construct:Monster} is based on a fixed
2B involution $x_{-1}$.  The centralizer of $x_{-1}$ is a maximal 
subgroup  $G_{x0}$ of $\setM$ of structure $2_+^{1+24}.\mbox{Co}_1$.
We use the ATLAS notation for describing the structure of a group, 
see \cite{Atlas}. The normal subgroup $Q_{x0}$ of structure $2_+^{1+24}$ 
of $G_{x0}$ is an extraspecial 2 group of plus type. There is a 
natural homomorphism $\lambda$ from $Q_{x0}$ to
$\Lambda / 2 \Lambda$, i.e. to the  Leech lattice $\Lambda$
modulo~2. The kernel of  $\lambda$ is equal to the 
centre $\{1, x_{-1}\}$ of $Q_{x0}$. 
The factor group $\mbox{Co}_1$ of $G_{x0}$ is the
automorphism  group of  $\Lambda / 2 \Lambda$. The automorphism
group $\mbox{Co}_0$ of the Leech lattice   $\Lambda$ has
structure  $2.\mbox{Co}_1$ with  centre of order 2. That centre is
generated by the mapping $x \mapsto -x$, for $ x \in \setR^{24}$.
So the operation of $\mbox{Co}_1$ (and hence also of $G_{x0}$)
on $\Lambda$ is defined up to sign. 

The (unique) minimal faithful real representation of the Monster
$\mathbb{M}$ has dimension 196883. We call that representation
$196883_x$. The 196884-dimensional rational representation $\rho$
of $\mathbb{M}$ constructed in \cite{Conway:Construct:Monster} has
matrix entries in $\mathbb{Z}[\frac{1}{2}]$; and it is a direct sum
of $196883_x$ and the trivial representation $1_x$. As a 
representation of $G_{x0}$ the representation $\rho$ splits as 
follows:
\begin{equation}
  \label{eqn:decom:rho}
     1_x \oplus 196883_x =
     \rho = 300_x \oplus 98280_x \oplus (4096_x \otimes  24_x),
\end{equation}
where the numbers in the names of the representations indicate
their dimensions.
Here $24_x$ is the natural 24-dimensional representation of the
automorphism group $\mbox{Co}_0$ of the Leech lattice, and $300_x$
is the symmetric tensor square of $24_x$. Representations and  
$4096_x$ and $98280_x$ will be described in 
Sections~\ref{ssect:maximal:G0} and~\ref{sect:rep:M}.
Ignoring signs, the basis vectors of monomial representation
$98280_x$ of $G_{x0}$ are in a one-to-one correspondence with 
the $2 \cdot 98280$ shortest nonzero vectors of the Leech 
lattice $\Lambda$.

In this paper we focus on the part $300_x$ of $\rho$; and to some
extent we also use $98280_x$. Since $300_x$ is the symmetric
tensor square of $24_x$, a vector in $300_x$ has a natural
interpretation as a symmetric matrix acting on the Euclidean
space $\setR^{24}$ spanned by the Leech lattice. In the sequel 
we identify $300_x$ with the space of these symmetric matrices.
The trivial part $1_x$ of representation $\rho$
is the subspace of $300_x$ spanned by the unit matrix $1_\rho$.

A four-group $\{1, x_{-1},  x_{\Omega},  x_{-\Omega} \}$ of the Monster
containing three commuting 2B involutions is also considered in
\cite{Conway:Construct:Monster}. The centralizer of that four-group
is a maximal subgroup $N_0$ of $\mathbb{M}$ of structure 
$2^{2+11+22}.(M_{24} \times \mbox{Sym}_3)$. The
factor $2^{2+11+22}$ describes a certain 2 group. The group $M_{24}$
is the Mathieu group acting as a permutation group on
a set $\tilde{\Omega}$ of 24 elements.  $\mbox{Sym}_3$ is the
symmetric permutation group of 3 elements. $\mbox{Sym}_3$ acts 
naturally on the set $\{ x_{-1},  x_{\Omega},  x_{-\Omega} \}$ 
of 2B involutions. The intersection $N_{x0} = N_0 \cap G_{x0}$ 
has index 3 in $N_0$ and is a maximal subgroup of  $G_{x0}$. 
$N_{x0}$ and acts monomially on $\rho$. A (considerably more 
detailed) diagram of the relevant subgroups of the Monster is
given in Figure~\ref{figure:subgroups:monster2}.

A set of generators of  $N_{x0}$ is given in
\cite{Conway:Construct:Monster}, together with another generator
$\tau \in N_0 \setminus N_{x0}$ of order 3 that cyclically 
exchanges $ x_{-1}$,  $x_{\Omega}$, and  $x_{-\Omega}$.
$\tau$ is called the {\em triality element}. The operation of all
these generators on $\rho$, and also sufficient information for 
effectively computing in  $N_0$, is given in
\cite{Conway:Construct:Monster}. For obtaining
a complete set of generators of $\mathbb{M}$ we just need another
element  $\xi$ of $G_{x0} \setminus N_{x0}$. Such a generator
$\xi$, and its operation on $\rho$, has been constructed in 
\cite{Seysen20}. Let $\Gamma$ be the set of all these generators
of $\setM$; and let $\Gamma^*$ be the set of all words in  $\Gamma$.

\subsection{Computing in the subgroup $G_{x0}$ of the Monster}
\label{ssect:Recognizing}

According to the {\em the Pacific Island model} in \cite{Wilson13}
we may assume that computations in $G_{x0}$ are easy (or at least
doable), while computations in $\mathbb{M}$ outside
$G_{x0}$ are difficult.

In Section~\ref{sect:Recogniz:Gx0} we will construct a vector
$v_1 \in \rho_{15}$ with the following properties:

\begin{itemize}
	\item
	The only element of $\setM$ fixing $v_1$ is the neutral element.
	\item 
	From $v_1 g$, with $g \in  \setM$ unknown, we can effectively check
	if $g$ is in $G_{x0}$ or not.
	\item 
	From $v_1 g$, with $g \in G_{x0}$ unknown, we can effectively
	compute a word $g'$ in $\Gamma^*$ with $g' = g$.
\end{itemize}

By our construction of $v_1$ in Section~\ref{sect:Recogniz:Gx0},
the first property of $v_1$ follows from a result 
in~\cite{LPWW-Computer_Monster}.
For testing our implementation this gives us the invaluable
advantage, that equality of two words $g_1, g_2$ in $\Gamma^*$ 
can be tested by comparing $v_1 g_1$ with $v_1 g_2$; here the 
correctness of this test follows from results that are independent
of the new algorithm in this paper.

 \subsection{The basic idea for computing in the Monster}
 \label{ssect:basic:idea}

Our goal is to recognize an unknown element $g$ of $\setM$ from 
the images $v_i g$ of a few vectors in $v_i$ in $\rho$. Here 
'recognizing' means finding a word $h$ in  $\Gamma^*$  that maps
each image $v_i g$ to its preimage $v_i$. If only the neutral 
element of $\mathbb{M}$ fixes all vectors $v_i$ then 
$g = h^{-1}$. So this leads to an algorithm for reducing
an arbitrary word of generators of a  $\mathbb{M}$ to a standard form,
and hence to an effective algorithm for the group operation 
in $\mathbb{M}$.

Here we give a brief overview of our reduction method based on this 
idea. More details will be given in Section~\ref{sect:strategy}.

We use a fixed pair $(\beta^+, \beta^-)$ of 2A involutions
in the subgroup $Q_{x0}$ of $\setM$ with 
$\beta^+ \beta^- = x_{-1}$. The centralizer of a 2A involution has 
structure $2.B$ (where $B$ is the Baby Monster group); and it fixes
a unique one-dimensional subspace of the representation $196883_x$. 
In \cite{Conway:Construct:Monster} a vector in $196883_x$ fixed by 
the centralizer of a 2A involution $t$ in $\mathbb{M}$ is
called  an {\em axis} of $t$. For any 2A involution 
$t \in \setM$ we will define a unique axis $\mbox{ax}(t)$ 
in the representation $\rho = 1_x \oplus 198883_x$.
Then $\setM$ is transitive on these axes; and the centralizer
of a 2A involution $t$ in $\setM$ is equal to the centralizer of 
$\mbox{ax}(t)$ in $\setM$. 
We put $v^+ = \mbox{ax}(\beta^+)$,  $v^- = \mbox{ax}(\beta^-)$;
and we write $H^+$ for the centralizer of $x_\beta$ (or of $v^+$)
in $\setM$.

The key idea of the new algorithm is to track images of the axes
$v^+, v^-$ under the action of an element of the Monster.
Given an image $v^+g$ of $v^+$ (with $g$ unknown) we present an 
effective algorithm for computing a word $h_1$ in $\Gamma^*$ with
$v^+g h_1 = v^+$. So we have $g h_1 \in H^+$.
This means that computation in the Monster can be reduced to 
computation in the group $H^+$ of structure $2.B$, which was
not possible before. 

In the next step, we deal with $g' = g h_1$.
Given an image $v^-g'$ of $v^-$ (with $g' \in H^+$ unknown) we 
present an effective algorithm for computing a word $h_2 \in H^+$ 
(given as a word in $\Gamma^*$) with $v^-g' h_2 = v^-$. 
Since $H^+$ fixes $h_2$, we have $v^+ g h_1 h_2 = v^+$ and 
$v^- g h_1 h_2 = v^-$. So $g h_1 h_2$ centralizes
both, $x_\beta$ and $x_{-\beta}$, and hence also the product
$x_{-1} =  x_\beta x_{-\beta}$.

So we have  $g h_1 h_2 \in G_{x0}$. Using the idea in
Section~\ref{ssect:Recognizing}, we may compute a $h_3 \in G_{x0}$
as a word in $\Gamma^*$ with $g h_1 h_2 h_3 =1$, provided that
we know $v_1 g$ (and hence also $v_1 g h_1 h_2$).

So given a triple $(v_1 g, v^+g, v^- g)$ of vectors in $\rho$,
we can effectively compute the element $g$ from that triple of 
vectors as a word in $\Gamma^+$. As we shall see in 
Section~\ref{sect:Recogniz:Gx0} et seq.,
it suffices if $(v_1 g, v^+g, v^- g)$ is given as a 
triple of vectors in $\rho_{15}$.

\subsection{Reducing an axis and the geometry of the Leech lattice}
 \label{ssect:Reducing:geometry}

In this subsection we discuss the basic geometric idea used for
finding an element of $\mathbb{M}$ that transforms an arbitrary
image of axis $v^+$ to $v^+$. Here we generously assume that 
the relevant geometric computations in the Leech lattice are 
doable. Details will be given in Section~\ref{sect:Orbits:Gx0}.

According to \cite{Norton98} there are twelve orbits of $G_{x0}$
on 2A axes. We may find representatives of all these orbits by 
multiplying axis $v^+$ with random elements of  $\setM$.
This way we find out that the projections of all these
representatives onto  $300_x$  are nonzero positive semidefinite
symmetric matrices. The elements of  $G_{x0}$  perform 
orthogonal transformations of the  matrices in $300_x$. Hence the 
projections of all axes onto $300_x$ are nonzero and positive
semidefinite. The multiset of the eigenvalues of such a  matrix
(counting each eigenvalue with the dimension of its eigenspace)
is determined by the orbit of $G_{x0}$ on the corresponding axis. 
It turns out that such a multiset of eigenvalues also 
characterizes the orbit of $G_{x0}$ on an axis.

A positive semidefinite matrix $A \in 300_x$ can be visualized as 
the (possibly degenerated) ellipsoid
$\{ x \in \setR^{24} \mid x A x^\top \leq 1\}$ in the space
spanned by the Leech lattice. The
shape of such an ellipsoid is determined by the eigenvalues of $A$;
so it is a property of the orbit of $G_{x0}$ on $A$.

For any axis we may analyse the eigenspaces of its
projection on $300_x$, leading to a wealth of geometric 
information related to the Leech lattice $\Lambda$, 
on which $G_{x0}$ operates naturally up to sign.
It turns out that these eigenspaces are (usually) spanned by 
rather short vectors of the Leech lattice. This means that we
have pretty good control over the geometry of an  axis when 
operating inside $G_{x0}$. 

For mapping an axis to another orbit of $G_{x0}$, we have to 
apply a power $\tau^{\pm 1}$ of the triality element
$\tau$ to that axis. Note that this operation may change the 
shape of the ellipsoid corresponding to the projection
of the axis onto $300_x$. 

In Section~\ref{sect:Orbits:Gx0} we will see that for a axis $v$  
the set $\{v \tau^k G_{x0} \mid k = \pm1\}$ of 
orbits of $G_{x0}$ depends on the  orbit $v N_{x0}$ of $v$ only.
In Section~\ref{ssect:maximal:G0} we will see that $N_{x0}$ can 
also be considered as the  centralizer of the standard
co-ordinate frame  of the Leech lattice. The images of that
frame correspond to the vectors of minimal norm 8 in 
$\Lambda / 2 \Lambda$, i.e. in the Leech lattice mod 2.

Thus for an axis $v$ the shapes of the ellipsoids 
corresponding to $v \tau^{\pm 1}$ are determined by the 
position of the ellipsoid corresponding to $v$ relative to 
the standard co-ordinate frame of the Leech lattice. Determining
these shapes in all relevant cases is far from trivial; but at
least we have a geometric idea how to proceed.

We say that two  axes have distance $l$ if any word in
$\Gamma^*$ transforming one axis into the
other axis contains at least $l$ powers of $\tau$. So, 
geometrically, our task is to rotate an ellipsoid in the
Leech lattice corresponding to an axis into a 'good'
position relative to the standard co-ordinate frame by applying
a transformation in $G_{x0}$. Here a position of an ellipsoid is
'good' if applying one of the transformations $\tau^{\pm 1}$
decreases the distance of the corresponding axis to the 
standard axis $v^+$. This way we may decrease the distance
of a given axis from the standard axis $v^+$ by a sequence
of operations in  $G_{x0}$ and in $\{\tau^{\pm 1}\}$. 
Finally, we obtain an axis that has distance 0 from axis
$v^+$. Then  an easy computation inside the group
$G_{x0}$ will map that axis to $v^+$.

So we have to find 'good' co-ordinate frames in $\Lambda$ (or 
vectors of minimal norm 8 in $\Lambda/ 2 \Lambda$), when an 
ellipsoid over the Leech lattice or, equivalently, a symmetric 
$24 \times 24$ matrix  $A \in 300_x$ is given. 
From a computational point of view we prefer to search
for such 'good' vectors in  $\Lambda/ 2 \Lambda$ with respect to a
given matrix $A \in 300_x$. It turns out that the projection of an
axis onto the subspace $98280_x$ of $\rho$ gives us useful 
information about the $98280$ shortest nonzero vectors in 
$\Lambda/ 2 \Lambda$, which we will also use for finding 'good' 
vectors in $\Lambda/ 2 \Lambda$.

\subsection{Dealing with a pair of axes}
 \label{ssect:Dealing:pair}

We call an axis {\em feasible} if it is equal to an image 
$v^- h$ of $v^-$ for some $h \in H^+$. 
So we have to find an element 
of $H^+$ that transforms a feasible axis to the axis 
$v^-$. For performing this task we use essentially the same
methods as in the previous subsection.

Put $H = H^+ \cap G_{x0}$.  By our construction of $\tau$ and $H^+$,
the group $H^+$ is generated by $H$ and $\tau$. 
Since $H$ is a subgroup of $G_{x0}$, we'll have good control over the
operation of $H$ on feasible axes. So we'll apply  $\tau^{\pm 1}$ 
to a feasible axis whenever we may decrease the distance between a 
feasible axis and $v^-$, in a similar way as in the previous 
subsection. Here it is important to note that there are not too many
orbits of $H$ on feasible axes. From \cite{mueller_2008} we will 
conclude that there are 10 such orbits.

\subsection{Implementation}
 \label{ssect:Implementation}

The implementation~\cite{mmgroup2020} contains highly optimized
C programs for dealing with the structures involved in the
construction of the representation $\rho$ in
\cite{Conway:Construct:Monster}. These structures include
the  binary Golay code and its cocode,  the Parker
loop, the Leech lattice $\Lambda$ (modulo 2 and 3), and also
the automorphism groups of these structures, as discussed in 
Sections \ref{sect:N0} and \ref{sect:Leech}.

It suffices if all computations in representation $\rho$
described in this paper are done modulo 15, i.e. in  $\rho_{15}$.
In  \cite{mmgroup2020} we also have highly 
optimized functions for  multiplying a vector in
$\rho_3$ or $\rho_{15}$ with a generator in $\Gamma$. (We are 
are a bit sloppy here, calling elements of $\rho_{15}$ also vectors,
although 15 is composite.) Multiplication of a vector in $\rho_{15}$
with any generator in $\Gamma$ costs less than 160 microseconds 
on the author's computer.

A reference implementation in Python for demonstrating the 
new reduction algorithm is presented in \cite{mmgroup_doc},
Section {\em Demonstration code for the reduction algorithm}.
The main function 
{\tt reduce\textunderscore monster\textunderscore element}
in that implementation reduces an element of the Monster.

%%%%%%%%%%%%%%%%%%%%%%%%%%%%%%%%%%%%%%%%%%%%%%%%%%%%%%%%%%%%%%%%%%%%%%%%%%%%%%
\section{The maximal subgroup $N_0$ of the Monster $\mathbb{M}$}
\label{sect:N0}
%%%%%%%%%%%%%%%%%%%%%%%%%%%%%%%%%%%%%%%%%%%%%%%%%%%%%%%%%%%%%%%%%%%%%%%%%%%%%%

The Monster has a maximal 2-local subgroup $N_0$ of structure 
$2^{2+11+22}.(M_{24} \times \mbox{Sym}_3)$ that is used in
the construction \cite{Conway:Construct:Monster}. We briefly recap
the structures given in \cite{Conway:Construct:Monster} that are 
required for understanding the generators and relations defining 
$N_0$. For  background we  refer to 
\cite{Aschbacher-Sporadic, Conway-SPLG, citeulike:Monster:Majorana}.

\subsection{The Golay code and its cocode}
\label{ssect:Golay}

Let $\tilde{\Omega}$ be a set of size~24 and construct the vector
space $\setF_2^{24}$ as $\prod_{i\in \tilde{\Omega}} \setF_2$.
A {\em Golay code} $\mathcal{C}$   is a 12-dimensional linear 
subspace of  $\setF_2^{24}$ whose smallest weight is $8$.
This characterizes the Golay code up to permutation. Golay code 
words have weight 0, 8, 12, 16, or 24. Code words of weight~8 
and~12 are called {\em octads} and {\em dodecads}, respectively.
The automorphism group of $\mathcal{C}$ is the Mathieu group 
$M_{24}$, which is quintuply transitive on the set $\tilde{\Omega}$.

We identify  the power set of $\tilde{\Omega}$ with
$\setF_2^{24}$ by mapping each subset of $\tilde{\Omega}$ to its
characteristic function, which is a vector in $\setF_2^{24}$,
as in \cite{Conway:Construct:Monster}.
So we may write $\tilde{\Omega}$ for the Golay code word 
containing 24 ones. For  elements $d, e$ of $\setF_2^{24}$ 
we write $d \cup e$, $d \cap e$ for  their union and intersection,
and $d+e$ for their symmetric difference.
We write $|d|$ for the cardinality of a set $d$.

We use the specific Golay code constructed in
\cite{Conway-SPLG}, Ch.~11,
which is also used in \cite{Seysen20}. The implementation
\cite{mmgroup2020} fixes a basis of that Golay code for computational
purposes; and it numbers the elements of the  set  $\tilde{\Omega}$ 
from 0 to 23.  Occasionally we write a subset of $\tilde{\Omega}$ 
representing an element of  $\mathbb{F}_2^{24}$ as a set of 
integers $0 \leq i < 24$, with the obvious meaning.
 
Let $\mathcal{C}^*$ be the cocode of $\mathcal{C}$, with scalar
product $\left<d, \delta \right>  \in  \mathbb{F}_2$ for
$d \in \mathcal{C}, \delta \in \mathcal{C}^*$. An element of 
$\mathcal{C}^*$ either has a unique representative of weight less 
than 4 in $\mathbb{F}_2^{24}$, or a set of six disjoint representatives 
of weight 4. Such a set of six representatives is called a 
{\em sextet}; the representatives in a sextet are called
{\em tetrads}.  For $\delta \in  \mathcal{C}^*$ let
$|\delta|$ be the minimum weight of  $\delta$; so we have
$0 \leq |\delta| \leq 4$.

\subsection{The Parker loop}
\label{ssect:Parker}

The Parker loop  $\mathcal{P}$ is a non-associative Moufang loop which
is a double cover of the Golay code written multiplicatively.
For any $d \in \mathcal{P}$ we write $\tilde{d}$ for the image of
$d$ in the Golay code  $\mathcal{C}$ as in \cite{Conway:Construct:Monster}
We  write 1 for the neutral element in  $\mathcal{P}$ and $-1$ for the
other preimage of the zero element of $\mathcal{C}$ in $\mathcal{P}$.
Set-theoretic operations
on elements of $\mathcal{P}$ are interpreted as operations on subsets of 
$\mathbb{F}_2^{24}$  corresponding to the images of these elements
in $\mathcal{C}$, as in \cite{Conway:Construct:Monster}.
So the intersection of two elements of $\mathcal{P}$ has a natural
interpretation as an element of $\mathbb{F}_2^{24}$ or of $\mathcal{C}^*$.
For  $i \in \tilde{\Omega}$ we abbreviate the cocode word $\{i\}$ of 
weight 1 to $i$. For $d, e, f \in  \mathcal{P}$ we have
\[
   d^2 = (-1)^{|d|/4} \, , \; \,
   d \cdot e = (-1)^{|d \cap e|/2} \cdot  e \cdot d , \; \,
   ( d \cdot e) \cdot f =  (-1)^{|d \cap e \cap f|}
    \cdot d \cdot (e \cdot f) \, .
\]

For practical computations an element  $d_i$ of the Parker loop
$\mathcal{P}$ is represented as a pair 
$(\tilde{d}_i, \mu_i)$ with $\tilde{d}_i \in \mathcal{C}$,  
$\mu_i \in \mathbb{F}_2$.  The implementation \cite{mmgroup2020}
defines  a {\em cocycle} 
$\theta: \mathcal{C} \times \mathcal{C} \rightarrow  \mathbb{F}_2$,
so that the product in  $\mathcal{P}$ is given by:
\[
(\tilde{d}_1, \mu_1) \cdot (\tilde{d}_2, \mu_2) = 
(\tilde{d}_1 + \tilde{d}_2, 
\mu_1 + \mu_2 + \theta(\tilde{d}_1, \tilde{d}_2) ) \;.
\]
Cocycles in a loop like  $\mathcal{P}$ are discussed in 
\cite{Aschbacher-Sporadic}, Chapter 4.
The cocycle $\theta$ satisfies the conditions in \cite{Seysen20}, 
Section~3.3, so that we may use the results in \cite{Seysen20}
for computations.
$\theta$ is quadratic in the first and linear in the second
argument; so it can also be considered as a function from 
$\mathcal{C}$ to $\mathcal{C}^*$. An element
$(\tilde{d}_i, \mu_i)$ of $\mathcal{P}$ is called {\em positive}
if $\mu_i = 0$ and  {\em negative} otherwise. We write $\Omega$ 
for the positive preimage of the Golay code word $\tilde{\Omega}$
in  $\mathcal{P}$. 
The centre of $\mathcal{P}$ is $\{ \pm 1, \pm \Omega\}$.

Let $\AutStP$ be the set of {\em standard automorphisms} of $\mathcal{P}$,
i.e. the set of automorphisms that map to an automorphism of  $\mathcal{C}$
in $M_{24}$ when factoring out $\{\pm1\}$. Any $\delta \in \mathcal{C}^*$
acts as a standard automorphism on  $\mathcal{P}$ given by
$d \mapsto (-1)^{\langle d, \delta \rangle} d$ for $d \in \mathcal{P}$; 
we call $\delta$ a {\em diagonal automorphism} of  $\mathcal{P}$.
An element $\pi$ of $\AutStP$ is called {\em even} if it fixes $\Omega$ 
and {\em odd} if it negates $\Omega$; for 
$\delta \in \mathcal{C}^* \subset \AutStP$ this agrees with the parity
of $\delta$ in  $\mathcal{C}^*$. The group
$\AutStP$ has structure $2^{12}.M_{24}$; the extension does not
split, and its normal subgroup of structure $2^{12}$ is isomorphic to
the group  $\mathcal{C}^*$ of diagonal automorphisms.

We follow the conventions in \cite{Conway:Construct:Monster} 
for denoting elements of $\mathcal{P}$, $\mathcal{C}$
and  $\mathcal{C}^*$:
\[
\begin{array}{ll}
	d,e,f   & \mbox{denote elements of $\mathcal{P}$
		or, loosely, of its homomorphic image $\mathcal{C}$}, \\
	\delta, \epsilon,  &
	\mbox{denote elements of $\mathcal{C}^*$},\\  
	i,j   & \mbox{denote elements of $\tilde{\Omega}$,
		also considered  as elements of  $\mathcal{C}^*$ 
		of weight 1}, \\
	d \cap e & 
	\mbox{denotes the subset 
		$ \; d \cap e \mbox{ of } \mathbb{F}_2^{24} \,,\,$ 
		usually considered  as an element of $\mathcal{C}^*$}, \\
	\pi, \pi', \pi'' & \mbox{denote elements of}  \AutStP . %,
	%\mbox{or, loosely, of its homomorphic image } M_{24}.         
\end{array}
\]

\subsection{The subgroups $N_0$ and $N_{x0}$ of the Monster}
\label{ssect:N0}

In   \cite{Conway:Construct:Monster} a subgroup $N_{x0}$ of $\mathbb{M}$ 
of structure  $2_+^{1+24}.2^{11}.M_{24}$ is defined. This group has
generators $x_\delta, x_d, y_d, x_\pi$, $\delta \in \mathcal{C}^*$, 
$d \in \mathcal{P}$, $\pi \in \AutStP$. In this paper we use the
generators of $N_{x0}$ defined in \cite{Seysen20}, with the same 
names as in \cite{Conway:Construct:Monster}, but with sightly different 
sign conventions, leading to simpler relations in $N_{x0}$. Generators
$x_d$ and $y_d$ in \cite{Seysen20} are equal to $x_d x_{-1}^{|d|/4}$ 
and $y_d y_{-1}^{|d|/4}$ in \cite{Conway:Construct:Monster},
respectively; the other generators are the same in in both papers. 
This simplification has been proposed in 
\cite{citeulike:Monster:Majorana}, Ch. 2.7.

\begin{Theorem}
\label{Thm:Nx0}	
In  $N_{x0}$ we have the following relations:
\begin{equation*}
\begin{array}{clll}
	& x_d x_e = x_{de} x_{d \cap e},
	\hphantom{mmmm}
	 &   y_d y_e = y_{de} x_{d \cap e},
	\hphantom{mmmm}
	& x_{\delta} x_{\epsilon} = 	x_{\delta \epsilon}, \\
	&  [x_d, y_e] = x_{d \cap e} x_{\Omega}^{|d \cap e|/2},  \quad
	& [x_d, x_\delta] = x_{-1}^{\langle d, \delta\rangle}, 	
	& y_d x_\delta = x_\delta Z_d
	  Z_{-1^{\vphantom t}}^{\langle d, \delta \rangle},\\
	& y_\Omega = x_{-1}, & y_{-1} = x_{-\Omega}, \\
	&  x_{\pi'} x_{\pi'\!'} =   x_{\pi' \pi'\!'}	, 	
	& x_d x_{\pi} = x_{\pi} x_{d^{\pi}} ,
	& y_d x_\pi = y_\pi x_{d^\pi} ,
\end{array}
\end{equation*}
for $d, e \in \mathcal{P}; \, \delta, \epsilon \in 
\mathcal{C}^* \subset  \AutStP;  \, \pi, \pi', \pi'' \in \AutStP$,
$\pi$ even. We have to put $Z_{\pm d}  = y_{\pm d}$ if $\delta$ 
is even, and $Z_{\pm d}  = (x_{\pm d} y_{\pm d})^{-1}$ if $\delta$  
is odd.  
\end{Theorem}

Theorem \ref{Thm:Nx0} follows from Theorem~1 and the definitions in \S 6 in
\cite{Conway:Construct:Monster}, or from Theorem~5.1 in in \cite{Seysen20}.
Here $[a,b]$ is the commutator $a^{-1}b^{-1} a b$; and  $d \cap e$ is an 
element of  $\mathcal{C}^*$.

We obtain a larger
group $N_0$ of structure $2^{2+11+22}.(M_{24} \times \mbox{Sym}_3)$ with 
$N_0:N_{x0} = 3$ by adding another generator $\tau$ called the
{\em triality element} in \cite{Conway:Construct:Monster}.
$\tau$ satisfies the relations:
\begin{equation}
	\label{eqn:tau}
   \tau^3 = 1,  \;  x_d \tau = \tau y_d, \;  y_d \tau = \tau (x_d y_d)^{-1}, \;
    [\tau, x_\pi] = 1,  \; x_{\pi'} \tau =  \tau^{-1} x_{\pi'}; \;
     \pi \mbox{ even},  \;
   \pi' \mbox{ odd} \; .
\end{equation}

The group  $N_0$  is a maximal subgroup of the Monster  $\mathbb{M}$.
The relations (\ref{eqn:tau}) can be obtained from the discussion in
\cite{Conway:Construct:Monster}, \S 6; or from \cite{Seysen20}, 
Theorem~5.1.
Computation in the group  $N_0$ is easy using the generators and 
relations given above. On the author's computer the group operation 
in  $N_0$ costs about 2 microseconds.

%%%%%%%%%%%%%%%%%%%%%%%%%%%%%%%%%%%%%%%%%%%%%%%%%%%%%%%%%%%%%%%%%%%%%%%%%%%%%%
\section{The Leech lattice and the maximal subgroup $G_{x0}$ of $\mathbb{M}$}
\label{sect:Leech}
%%%%%%%%%%%%%%%%%%%%%%%%%%%%%%%%%%%%%%%%%%%%%%%%%%%%%%%%%%%%%%%%%%%%%%%%%%%%%%

\newcommand{\InvSqrtEight}{ 
	{ \frac{\scriptstyle 1}{ \sqrt{\scriptstyle 8}}}}

\newcommand{\SmallTextMbox}[1]{\mbox{\rm{\scriptsize #1}}}

\subsection{The Leech lattice and its relation to the subgroup $Q_{x0}$ of $N_{x0}$}
\label{ssect:Leech:Qx0}

Let  $\mathcal{C}$ be the Golay code in $\Field_2^{24}$ as in
Section \ref{ssect:Golay}.
Let  $\{\eta_i,  i \in \tilde{\Omega}\}$ be a basis
of the Euclidean vector space $\setR^{24}$ so that the basis vectors 
of both, $\Field_2^{24}$ and $\setR^{24}$, are labelled by the elements 
of the  same set $\tilde{\Omega}$. Then the Leech lattice $\Lambda$ is 
the set of vectors 
$u = \sum_{i \in \tilde{\Omega}} u_i \eta_i$, $u_i \in \setZ$,  such
that there is an $m \in \{0,1\}$ and a $d \in \mathcal{C}$ with 
\[
\begin{array}{rllcl}
	\forall \, i \in {\tilde{\Omega}}  : & u_i
	&= &m + 2\cdot \left< d, i \right> &\pmod{4} \; ,\\
	\textstyle \sum_{i \in {\tilde{\Omega}}} & u_i
	&=&  4m &\pmod{8}  \; . 
\end{array}
\]
Here we scale the basis vectors $\eta_i$ of $\setR^{24}$ so that 
they have length $\InvSqrtEight$ and not 1. Then $\Lambda$ is 
the unique even unimodular lattice in $\setR^{24}$
such that the shortest nonzero vectors have squared norm 4, 
see e.g. \cite{Conway-SPLG}, Ch.~4.11 for background.
Thus for vectors $u, v \in \Lambda$ with co-ordinates $u_i, v_i$
the scalar product  $\left< u, v\right>$ is equal to  
$\frac{1}{8}\sum_{i \in \tilde{\Omega}} u_i v_i$. For $u \in \Lambda$
we define $\mbox{type}(u)$ = $\frac{1}{2}\left<u, u\right>$; 
so a shortest nonzero vector in $\Lambda$ is of type~2.

The group  $N_{x0}$ has a normal subgroup  $Q_{x0}$ of
structure $2_+^{1+24}$, generated by
$x_d, x_\delta$, $d \in \mathcal{P}, \delta \in \mathcal{C}^*$.
So $Q_{x0}$ is an extraspecial 2 group of order  $2^{1+24}$ of plus
type. The centre of $Q_{x0}$  is  $\{x_{\pm 1}\}$.
Let $\Lambda/2\Lambda$ be the Leech lattice modulo 2.
There is a homomorphism $\lambda$ from $Q_{x0}$ onto
$\Lambda/2 \Lambda$, with kernel $\{x_{\pm 1}\}$, given by:
\begin{equation}
\label{eqn:lambda:Leech}
    \textstyle
    x_d \stackrel{\lambda}{\longmapsto}
     {\textstyle \frac{1}{2} \sum_{j \in d} \lambda_j}  , \; \;
    x_i \stackrel{\lambda}{\longmapsto}  \lambda_i ,  \quad
     \mbox{for} \;  d \in \mathcal{P}, i \in \tilde{\Omega},
     \quad \mbox{where} \;
    \lambda_i = -4 \eta_i + \sum_{j \in \tilde{\Omega}} \eta_j \, .	  
\end{equation}
We write
\[
\begin{array}{ll}
	 x_r, x_s & \mbox{ for general elements of } Q_{x0} \, , \\
	 \lambda_r, \lambda_s &  \mbox{ for the elements  }
	   \lambda({x_r}), \lambda({x_s})  \mbox{ of } \Lambda / 2 \Lambda \, .
\end{array}
\]
For $x_r, x_s \in Q_{x0}$ we have
\[
[x_r,x_s]  =  x_{-1}^{\left<\lambda_r, \lambda_s\right>} \; ,  \quad
     x_r^2 = x_{-1}^{\SmallTextMbox{type}(\lambda_r)}.
\]
The mapping $\lambda_r \mapsto \mbox{type}(\lambda_r) \pmod{2}$ defines
the natural non-singular quadratic form on  $\Lambda/ 2\Lambda$. For a 
proof of these statements see \cite{Conway:Construct:Monster}, 
Theorem~2, or \cite{Seysen20}, Theorem~6.1. 

The type of a vector in $\Lambda/ 2\Lambda$ is the type of its shortest
preimage in  $\Lambda$. Every vector in  $\Lambda/ 2\Lambda$ 
has type 0, 2, 3, or 4. Vectors of type 2 are called {\em short};
there are 98280 short vectors in  $\Lambda/ 2\Lambda$.
The automorphism groups of  $\Lambda$ and  $\Lambda/ 2\Lambda$ are
called $\mbox{Co}_0$ and  $\mbox{Co}_1$. The group  $\mbox{Co}_1$
is simple; and  $\mbox{Co}_0$ has structure  $2.\mbox{Co}_1$.
$\mbox{Co}_0$ and $\mbox{Co}_1$ are transitive on the sets of
vectors of types 2, 3, or 4. A vector of type 2 or 3 in
$\Lambda/ 2\Lambda$ has two opposite preimages of the same
type in $\Lambda$. 

The vector $\lambda_\Omega=\lambda(x_\Omega)$ in
$\Lambda/ 2\Lambda$ is of type 4; its preimages of type~4 in $\Lambda$
are the 48 vectors of type 4 proportional to the unit vectors 
in $\setR^{24}$. The images of the standard co-ordinate frame 
of  $\Lambda$ under $\mbox{Co}_0$  are in one-to-one 
correspondence with the vectors of type~4
in  $\Lambda/ 2 \Lambda$, see  \cite{Conway-SPLG}, Ch. 10.3.3.

For  $\lambda_r \in \Lambda / 2 \Lambda$ put $\Lambda(\lambda_r) = 
\{ v \in \Lambda \mid v = \lambda_r \pmod{ 2 \Lambda}, \mbox{type}(v) 
=  \mbox{type}(\lambda_r)\}$. So  $\Lambda(\lambda_r)$ is the set of
the shortest preimages of  $\lambda_r$ in $\Lambda$.
E.g. $\Lambda(\lambda_\Omega)$ is the set 
$\{\pm 8 \eta_i \mid i \in \tilde{\Omega}\}$.

\subsection{The maximal subgroup $G_{x0}$ of the Monster}
\label{ssect:maximal:G0}

We construct a maximal subgroup $G_{x0}$ of structure
$2_+^{1+24}.\mbox{Co}_1$ of $\setM$ with $G_{x0} \cap N_0 = N_{x0}$
as follows. The extraspecial normal subgroup $Q_{x0}$ of $N_{x0}$ 
has a unique irreducible real representation $4096_x$ of dimension 4096,
see \cite{Conway:Construct:Monster, citeulike:Monster:Majorana}.
That representation may be extended to a representation of a (unique)
group of structure $2_+^{1+24}.\mbox{O}_{24}^+(2)$, where 
$\mbox{O}_{24}^+(2)$ is the orthogonal group on  $\setF_2^{24}$ of plus 
type. Since  $\mbox{Co}_1 \subset \mbox{O}_{24}^+(2)$, we also obtain 
a representation  $4096_x$ of a group  $G'_{x0}$ of structure  
$2_+^{1+24}.\mbox{Co}_1$. 

If $G_1, G_2$ are groups with a common factor group $H$
and homomorphisms $\phi_i: G_i \rightarrow H$, $i=1,2$, 
then the {\em fibre product} $G_1 \bigtriangleup_{H} G_2$ is
the subgroup of the direct product $G_1 \times G_2$
defined by:
\[
G_1  \bigtriangleup_{H} G_2 \; = \;  
\left\{(x,y)   \in G_1 \times G_2  
\mid  \phi_1(x) = \phi_2(y)  \right\} \; .
\]  
We define $G_{x0}$  by 
\[
   G_{x0} = {\textstyle \frac{1}{2}}
  (  G'_{x0} \bigtriangleup_{\mbox{\scriptsize Co}_1} \mbox{Co}_0 ).
\]
Here the factor $\frac{1}{2}$ means that we identify the centres 
(of order 2) of the groups $G'_{x0}$ and  $\mbox{Co}_0 $.
We remark that $G_{x0}$ and $G'_{x0}$ are not isomorphic. 

From \cite{Conway:Construct:Monster} or from Theorem~\ref{Thm:Nx0}
we see that $N_{x0}$ normalizes the four-group 
$\{x_{\pm 1}, x_{\pm \Omega} \}$. By \cite{Atlas} the group
$2^{11}. M_{24}$ is maximal in $\mbox{Co}_1$; so we conclude 
that $N_{x0}$ is the normalizer of that four-group in $G_{x0}$, 
and hence also the stabilizer of 
$\lambda_\Omega = \lambda(x_{\pm \Omega}) \in \Lambda / 2\Lambda$
in $G_{x0}$. So we may say that  $N_{x0}$ is the stabilizer of 
the standard co-ordinate frame of  $\Lambda$.

Computation in $\mbox{Co}_0$ is easy. Computation in $4096_x$ can be 
greatly accelerated by using the fact that a certain basis of the 
underlying vector space $V = \setR^{4096}$ has a natural structure as a
12-dimensional affine space over $\mathbb{F}_2$. Using such a basis, 
the standard basis of $\hom(V, V)$ inherits the structure of 
24-dimensional affine space $A$ over $\mathbb{F}_2$. It turns out that 
a matrix in $\hom(V,V)$ representing a group element has
entries with a fixed absolute value in an affine subspace of $A$,
and entries 0 elsewhere.  The signs of the nonzero
entries of the matrix are essentially given by a quadratic form 
on that subspace. Using these ideas, the group operation in $G_{x0}$ 
can be done in a bit more than 10 microseconds on the author's computer.
%Computing of a character in  $4096_x$ is also easy.
For details we refer to the documentation of the project
\cite{mmgroup2020}.

Using the functions in \cite{mmgroup2020} we may compute the
character of any element of the subgroup $G_{x0}$ of $\mathbb{M}$ in 
the 196883-dimensional real representation of $\mathbb{M}$ in a few
ten milliseconds.

%%%%%%%%%%%%%%%%%%%%%%%%%%%%%%%%%%%%%%%%%%%%%%%%%%%%%%%%%%%%%%%%%%%%%%%%%%%%%%
\section{The 196884-dimensional representation of $\mathbb{M}$}
\label{sect:rep:M}
%%%%%%%%%%%%%%%%%%%%%%%%%%%%%%%%%%%%%%%%%%%%%%%%%%%%%%%%%%%%%%%%%%%%%%%%%%%%%%

Let $24_x$ be the representation of the group $\mbox{Co}_0$ on 
$\setR^{24}$ as the automorphism group of the Leech lattice $\Lambda$, 
and let $4096_x$ be as in Section~\ref{ssect:maximal:G0}.  Then the 
maximal subgroup $G_{x0}$ of $\mathbb{M}$ has a faithful real 
representation $4096_x \otimes 24_x$. A construction of 
$G_{x0}$  and of its representation 
$4096_x \otimes 24_x$ is given in \cite{Conway:Construct:Monster},
where the generators $x_\delta, x_d, y_d, x_\pi$ of the maximal
subgroup $N_{x0}$ of $G_{x0}$ are given explicitly.
In  \cite{Seysen20} we define another generator  $\xi$ of 
order 3 in $G_{x0} \setminus N_{x0}$, and its action on
$4096_x \otimes 24_x$. Using the basis of $\Lambda$ in 
$\setR^{24}$ in Section~\ref{ssect:Leech:Qx0}, the generator 
$\xi$ acts on $\setR^{24}$ by right multiplication 
with the matrix
\begin{equation}
\label{eqn:matrix:xi}	
   \renewcommand{\arraystretch}{0.85}
   \left(
   \begin{array}{cccccc}
   \scriptstyle \! \! \! \!  AB  \! \! \! \! & \\
   	     &\scriptstyle  \! \! \! \!   AB \! \! \! \!   \\
   	     & & \scriptstyle \! \! \! \!   AB  \! \! \! \!   \\
   	     & & & \scriptstyle \! \! \! \!   AB  \! \! \! \!   \\
   	     & & & &\scriptstyle  \! \! \! \!   AB  \! \! \! \!   \\
   	     & & & & &  \scriptstyle \! \! \! \!   AB  \! \! \! \!   \\
   \end{array}
   \right)  , 
    \; \;  \textstyle A = \frac{1}{2}
   \left(
   \begin{array}{cccc}
	 \! \! \scriptstyle -1 \! \! & \scriptstyle 1 & \scriptstyle 1 & \scriptstyle 1 \\
     \scriptstyle 1 & \! \! \scriptstyle -1 \! \! & \scriptstyle 1 & \scriptstyle 1 \\
     \scriptstyle 1 & \scriptstyle 1 & \scriptstyle \! \! -1 \! \! & \scriptstyle 1 \\
     \scriptstyle 1 & \scriptstyle 1 & \scriptstyle 1 & \scriptstyle \! \! -1 \! \!\\
   \end{array}
   \right)  , \;
     \textstyle B =
   \left(
   \begin{array}{cccc}
   \! \!  \!	\scriptstyle -1 \! \!  \! & \\
    &   \! \scriptstyle  \!   1  \\
	& &   \!\scriptstyle   1  \! \\
	& & &   \! \scriptstyle  1  \! \\
   \end{array}
   \right)  ,
   \renewcommand{\arraystretch}{1.0}
\end{equation}
up to sign.
A suitable basis of $4096_x \otimes  24_x$ is
given in \cite{Conway:Construct:Monster} and also in \cite{Seysen20},
with slightly different sign conventions leading to a simpler
construction of generator $\xi$. We omit the
details, since we do not need them in this paper.

The group $G_{x0}$ operates on its normal subgroup $Q_{x0}$ by 
conjugation.
% For $x_r \in Q_{x0}$ let 
% $\mbox{type}(x_r) = \mbox{type}(\lambda(x))$. Then for $q \in Q_{x0}$ we
% have $\mbox{type}(x) = \mbox{type}(x^q)$. 
Let $98280_x$ be a real vector space, with basis vectors 
$X_{d \cdot \delta}$,  $d \in \mathcal{P}$, 
$\delta \in  \mathcal{C}^*$, such that 
$\mbox{type}(\lambda(x_d \cdot x_\delta)) = 2$. We identify
$X_{-d \cdot \delta}$ with $-X_{d \cdot \delta}$. Then an
element of $G_{x0}$ operates on $X_{d \cdot \delta}$ in the same
way as it operates on $x_d \cdot x_\delta$ by conjugation.
With this operation the space  $98280_x$ is a
98280-dimensional monomial representation of $G_{x0}$
with kernel $\{1, x_{-1}\}$.

For constructing the 196884-dimensional representation $\rho$ of 
$\mathbb{M}$ we also need the symmetric tensor square
$300_x = 24_x \otimes_{\mbox{\scriptsize sym}} 24_x$. 
Note that $24_x$ is not a representation 
of $G_{x0}$, but $300_x$ is. Basis vectors of $24_x$ are
$\eta_i$, $i \in\tilde{\Omega} = \{0,\ldots,23\}$,
as in Section~\ref{ssect:Leech:Qx0}. Then 
an element of  $300_x$ has a natural interpretation as symmetric
real $24 \times 24$ matrix. For $i, j \in \tilde{\Omega}$ we write 
$\eta_{i,j}$ for the symmetric matrix with an entry 1 in row $i$, 
column $j$, and also in row $j$, column $i$, and zeros elsewhere.
So the elements of the set 
$\{ \eta_{i,j} \mid i, j \in \tilde{\Omega},  i \leq j\}$ 
form a basis of $300_x$.

We define the representation $\rho$ of $G_{x0}$ by:
\begin{equation}
 \label{eqn:rho}
    \rho = 300_x \oplus 98280_x \oplus (4096_x \otimes  24_x) \, .
\end{equation}

For extending $\rho$ to a representation of $\mathbb{M}$, Conway 
\cite{Conway:Construct:Monster} defines the action of the generator
$\tau \in \mathbb{M} \setminus G_{x0}$ of order~3 on $\rho$. 
For a proof that $\rho$ actually represents $\mathbb{M}$, he has 
to show that there is a certain algebra on $\rho$
invariant under $\mathbb{M}$. Therefore he constructs
an algebra visibly invariant under $G_{x0}$; and he shows that this
algebra is also invariant under $\tau$. That algebra is called the
{\em Griess algebra}; it has first been constructed by
Griess \cite{Griess:Friendly:Giant}.

Representation $300_x$ can also be decomposed as 
$300_x = 1_x \oplus 299_x$, where $1_x$ is the trivial representation
of $G_{x0}$ corresponding to the subspace of $300_x$ spanned by the 
unit matrix $1_\rho$. Then $299_x$ is an irreducible 
representation of  $G_{x0}$ containing the symmetric $24 \times 24$ 
matrices in $300_x$ with trace zero. Replacing $300_x$ by $299_x$ 
in (\ref{eqn:rho}) we obtain the smallest faithful irreducible 
representation $196883_x$ of  $\mathbb{M}$ in characteristic~0. 

Representation $\rho$ preserves a positive-definite quadratic form.
With respect to that form, basis vectors $\eta_{i,j}$ of $\rho$
have squared norm 2 for $i \neq j$, and all other basis vectors 
have norm 1. Any two basis vectors are orthogonal, unless 
equal or opposite.

The following set $\Gamma$ generates 
$N_{x0}$ and $G_{x0}$, and hence the Monster group $\mathbb{M}$:
\[
\Gamma = \left\{ x_\delta, x_d, y_e, x_\pi, \tau^{\pm 1}, 
 \xi^{\pm 1} \mid
\delta \in \mathcal{C}^*; \,  d, e \in \mathcal{P}; \, 
\pi \in \AutStP; 
\right\} \; .
\]
The operation of the generators in $\Gamma$ on $\rho$ (except 
for $\tau^{\pm 1}$ and $\xi^{\pm 1}$) is given in Table~1 in 
\cite{Conway:Construct:Monster}, or in Table~3 in \cite{Seysen20}. 
The operation of  $\tau^{\pm 1}$ is also given in Table~3 in 
\cite{Seysen20};
in \cite{Conway:Construct:Monster} this corresponds to a cyclic 
exchange of the three languages in the {\em dictionary} in Table~2.
The operation of $\xi$ on $\rho$ is described in \cite{Seysen20},
Section 9. We write $\Gamma^*$ for the set of words in $\Gamma$.

The matrices in  $\rho$ representing the generators 
in $\Gamma$ are sparse matrices. Any such matrix can be represented 
as a product of  monomial matrices, and a sequence of at most six 
matrices of Hadamard type. Here a matrix of Hadamard type is a matrix 
with diagonal blocks and blocks of $2 \times 2$ matrices of shape
$c \begin{psmallmatrix} 1 & 1\\ 1 & -1\end{psmallmatrix}$
for $c \in \{\frac{1}{2}, 1\}$. For any odd integer $p > 1$ we  
define the representation $\rho_p$ of the Monster as the 
representation $\rho$, with matrix entries taken modulo $p$. 
$\rho_p$ is well defiend for odd $p$, since the denominators of
the entries of the matrices in $\rho$  are powers of two. 
%In this paper we deal with $\rho_3$ and $\rho_{15}$.

%%%%%%%%%%%%%%%%%%%%%%%%%%%%%%%%%%%%%%%%%%%%%%%%%%%%%%%%%%%%%%%%%%%%%%%%%%%%%%
\section{Recognizing an element of $G_{x0}$ in  $\mathbb{M}$}
\label{sect:Recogniz:Gx0}
%%%%%%%%%%%%%%%%%%%%%%%%%%%%%%%%%%%%%%%%%%%%%%%%%%%%%%%%%%%%%%%%%%%%%%%%%%%%%%

In this section we will show that computing in the maximal
subgroup  $G_{x0}$ of  $\mathbb{M}$ is easy.

Given an element $g$ of $\mathbb{M}$ as a word in $\Gamma^*$, we 
want to check if $g$  is in  $G_{x0}$ or not. If this is the case 
then we also want to obtain a representation of $g$ as a word in the
generators in $\Gamma$ that are in $G_{x0}$. That representation 
should depend on the value of $g$ only, and not on the given 
representation of $g$ as a word in $\Gamma^*$. In this section 
we construct a $v_1 \in \rho_{15}$ such that every $g\in G_{x0}$ can 
effectively be reconstructed as a word in the generators of  $G_{x0}$ 
from $v_1 \cdot g$. The algorithm for reconstructing $g$ also 
detects if $g$ is in $G_{x0}$ or not. 

An algorithm for solving that problem must certainly be able to detect
if a word in $\Gamma^*$ is the identity or not. For this purpose we may
use two nonzero vectors $v_{71}$ and $v_{94}$  in the representation
$198883_x$ of $\mathbb{M}$ that are fixed by an element of order 71 
and negated by an element of order 94, respectively. 
In \cite{LPWW-Computer_Monster} it is shown that only the
neutral element of  $\mathbb{M}$  fixes both vectors $v_{71}$ and
$v_{94}$.
So we may check if an element is the identity in the Monster. Here the 
corresponding calculations can be done modulo any odd prime; and we 
may even use different primes for the operation on  $v_{71}$ and on
$v_{94}$. Since we actually compute in the representation $\rho$ 
of $\mathbb{M}$, we must make sure that the projections of the 
vectors $v_{71}$ and $v_{94}$ from $\rho$ onto the subspace
$196883_x$ are not zero.

We generate $v_{71}$ as a vector in the representation $\rho_3$ of the 
Monster in characteristic~3, and  $v_{94}$ as a vector in the 
representation $\rho_5$ of the Monster in characteristic~5. We combine 
these two vectors to a vector $v_1$ in the representation $\rho_{15}$ 
via Chinese remaindering. Note that a computation in $\rho_{15}$ is 
faster than the combination of two similar computations in $\rho_3$. 
We will impose additional restrictions onto the vector $v_3$, so
that membership in  $G_{x0}$  can be tested.

In Section \ref{sect:rep:M} we have seen that an element of 
the subspace $300_x$ of $\rho$ corresponds to a symmetric 
matrix on the space $\setR^{24}$ containing  $\Lambda$. 
For $v \in \rho$ we write $M(v)$ for the symmetric matrix on 
$\setR^{24}$ corresponding to the projection of $v$ onto $300_x$.
When searching for a vector $v_{71} \in \rho_3$ fixing an element of
order $71$, we impose an additional condition on the matrix 
$M(v_{71})$. Modulo 3, we require that $M(v_{71})$ has corank 1, and 
that the kernel of $M(v_{71})$ contains a basis vector of 
$\setR^{24}$. In Appendix~\ref{app:probability} we estimate the 
cost of finding a suitable vector $v_{71}$. On the author's 
computer this takes less than a minute in average.

Given a fixed suitable $v_{71} \in \rho_3$ and a $g \in \mathbb{M}$ 
as a word in  $\Gamma^*$, we can compute the kernel of 
$M(v_{71} \cdot  g)$  modulo 3. 
In case $g \in G_{x0}$ that kernel is spanned by the image of a 
basis vector of $\Lambda$ (modulo 3); and we can easily compute
the image of that basis vector in $\Lambda$ under $g$ (modulo~3)
from $M(v_{71} \cdot  g)$ up to sign. 
So we can also compute the image $\lambda' = \lambda_\Omega g$ 
of the standard co-ordinate frame $\lambda_\Omega$ in 
$\Lambda / 2\Lambda$. In Appendix~\ref{app:map:type4} we 
show how to compute a $h_1 \in G_{x0}$ with 
$\lambda' h_1  = \lambda_\Omega$, for any
$\lambda' \in \Lambda/2\Lambda$ of type~4. So $g h_1$ fixes 
$\lambda_\Omega$. Since $N_{x0}$ is the centralizer of the 
standard frame $\lambda_\Omega$ in $G_{x0}$, we obtain 
$g h_1 \in N_{x0}$.

So we are left with the identification of an element $g$ of
$N_{x0}$ from a vector $v_{1} g$, where  $v_{1} \in \rho_{15}$ 
is fixed as above. The factor group $M_{24}$ of $N_{x0}$ acts
on the matrix $M(v_1)$ by permuting the rows and
columns of  $M(v_1)$ as given by the natural permutation
action of $M_{24}$, up to sign changes in the matrix.
Thus $M_{24}$ permutes the multisets given by the absolute 
values of the entries of a row. For a random matrix 
$M(v_1 \cdot g), g \in N_{x0}$, (with entries given modulo 15)
we can almost certainly 
recover the permutation in $M_{24}$ corresponding to $g$
from the operation on these multisets. 
In the implementation we compute a hash function on the 24 
multisets given by the absolute values of the entries of the 
rows of $M(v_1)$; and we precompute a new vector $v_1$ if 
these 24 hash values are not mutually different. So we can easily 
find a generator $x_\pi$, $\pi \in \AutStP$, such that 
$g x_\pi^{-1}$ is in the kernel of structure $2^{1+24+11}$ of 
the homomorphism $N_{x0} \rightarrow M_{24}$.

An element $g$ in the group $2^{1+24+11}$ given above can be
reduced to an element of  $Q_{x0} = 2_+^{1+24}$ by a sequence 
of sign checks in the matrix  $M(v_1 \cdot g)$. In the 
unlikely case that  $M(v_1)$ has too many zero entries 
we just precompute another $v_1 \in \rho_{15}$.

Finally, an element of the extraspecial 2-group  $Q_{x0}$
can be recognized by a sequence of sign checks in the parts 
of  $v_1$  that are not contained $300_x$.  At the end we
arrive at a relation $g h = 1$ for a word $h$ in the
generators of $G_{x0}$; and we have to check  that $g h$
really fixes $v_1$. 

The algorithm sketched in this section succeeds if the
input $g$ is in $G_{x0}$ and fails otherwise. 
%In case of failure we abort as early as possibly in order to 
%minimize the average run time in that case.
For an implementation of that algorithm we refer to
function {\tt reduce\textunderscore G\textunderscore x0}
in \cite{mmgroup_doc},
Section {\em Demonstration code for the reduction algorithm}.
 
%%%%%%%%%%%%%%%%%%%%%%%%%%%%%%%%%%%%%%%%%%%%%%%%%%%%%%%%%%%%%%%%%%%%%%%%%%%%%%
\section{The strategy for reducing an element of $\mathbb{M}$}
\label{sect:strategy}
%%%%%%%%%%%%%%%%%%%%%%%%%%%%%%%%%%%%%%%%%%%%%%%%%%%%%%%%%%%%%%%%%%%%%%%%%%%%%%

In this section we explain the idea how to obtain an unknown
element $g$ of the Monster as a word in $\Gamma^+$  from a triple
$(v_1 g, v^+ g, v^- g)$, where $v_1, v^+$, and $v^-$ are fixed 
elements of $\rho_{15}$.  

\subsection{Two 2A involutions and their centralizers $H^+, H^-$ in $\setM$}
\label{ssect:central:H}

The Monster $\mathbb{M}$ has two classes of involutions, which 
are called 2A and 2B in the ATLAS \cite{Atlas}. 
%The group $G_{x0}$ is the centralizer of the 2B involution $x_{-1}$ 
%in $\mathbb{M}$. 
The centralizer of a 2A involution is a group of structure $2.B$,
where $B$ is the Baby Monster. 
An element $t$ of the subgroup $Q_{x0}$ of $\mathbb{M}$ is a
2A involution in $\setM$ if and only if $\lambda(t)$ is of
type~2 in $\Lambda / 2 \Lambda$, see
\cite{Aschbacher-Sporadic, Conway:Construct:Monster,
	citeulike:Monster:Majorana}.
In the remainder of this paper let $\beta$ be the fixed element of 
the Golay cocode $\mathcal{C}^*$ given by the subset $\{2,3\}$ 
of $\tilde{\Omega} = \{0,\ldots, {23}\}$, using the basis 
of $\mathcal{C}$ in \cite{mmgroup2020}. Then $x_\beta \in Q_{x0}$.
Put $x_{-\beta} = x_{-1} x_\beta$, where $x_{-1}$ is the central 
involution in $Q_{x0}$ as in Section~\ref{ssect:Leech:Qx0}.

Then $x_\beta$ and $x_{-\beta}$ are 2A involutions in  $\mathbb{M}$;
and $\lambda_\beta := \lambda(x_\beta) = \lambda(x_{-\beta}) $  
is of type~2 in $\Lambda / 2\Lambda$.
Let $H^+$ and $H^-$ be the centralizers of $x_\beta$
and $x_{-\beta}$, respectively; and let $H = H^+ \cap H^-$. Since 
$H$ centralizes  $x_{-1} = x_\beta \cdot x_{-\beta}$,  we have 
$H  \subset G_{x0}$. $x_\beta$ centralizes $\xi$ and $\tau$,
see \cite{Seysen20}.

One can  show that $H = H^+ \cap G_{x0}$, and that $H$
has structure  $2^{2+22}.\mbox{Co}_2$, where $\mbox{Co}_2$ is
the (simple) subgroup of $\mbox{Co}_1$ fixing the type-2 vector
$\lambda_\beta$ in $\Lambda / 2 \Lambda$. 

A few words for justifying our choice of $\beta \in \mathcal{C}^*$ 
are appropriate. We require that $x_\beta$ is a 2A involution in
$\mathbb{M}$. For simplifying computations in the centralizer $H^+$
of  $x_\beta$ we also want to have $\tau, \xi \in H^+$ for the
generators $\tau, \xi$. This restricts $\beta$ to an element of 
$\mathcal{C}^*$ of weight~2 which is {\em coloured} in the 
terminology of \cite{Seysen20}. The coloured cocode words 
correspond to the elements of the cocode of the {\em hexacode},
which is a 3-dimensional linear code in $\mathbb{F}_4^6$ with 
Hamming distance~4, as defined in \cite{Conway-SPLG}, Ch.~11.2. 
Here our Golay cocode word $\beta$ corresponds to the element of
the cocode of the hexacode given by 
$(1,0,\ldots,0)  \in \mathbb{F}_4^6$, which is a fairly
natural choice.

\begin{figure}[!h]
	\centering
	\begin{tikzpicture}
		\matrix (m) [matrix of math nodes,row sep=1.75ex,column sep=0.2em,minimum width=2em]
		{  
			&          &  \mathbb{M} &           &         \\
			H^+  \! = 2.B &   \!\!\!\!\!\! H^- \!\! = \! 2.B   \!\!\!\!\!\! \\
			&    & G_{x0} = 2_+^{1+24}.\mbox{Co}_1	&  &
			\!  \!   N_0 \! = \! 2^{2+11+22}\!.(M_{24}\! \times \! \mbox{Sym}_3)  \\   
			H = 2^{2+22}\!.\mbox{Co}_2 \\
			&  &  N_{x0} \! = \! 2_+^{1+24} \! .2^{11} \! .M_{24}    \\     
			H \! \cap \! N_{x0} \! = \! 2^{2+22} \!.2^{10}\! .M_{22}.2 \\
			& & Q_{x0} = 2_+^{1+24}	&   &    \AutStP = 2^{12}.M_{24}\\    
			H \cap Q_{x0} = 2^{2+22} \\
			&  &    \mathcal{C}^* = 2^{12}  \\
			& \\
			& & 1 \\
		};         	
		\path[->] 
		(m-2-1) edge (m-1-3) 
		(m-3-3) edge  node[left] {$\scriptstyle \tau$} (m-1-3)
		(m-3-5) edge  node[above] {$\scriptstyle \xi$} (m-1-3)
		(m-4-1) edge  node[left] {$\scriptstyle \tau$}  (m-2-1)
		(m-4-1) edge  (m-3-3)
		(m-5-3) edge  node[left] {$\scriptstyle \xi$} (m-3-3)
		(m-5-3) edge  node[above] {$\scriptstyle \tau$}  (m-3-5)
		(m-6-1) edge  node[left] {$\scriptstyle \xi$} (m-4-1)
		(m-6-1) edge (m-5-3)
		(m-7-3) edge  node[left] {$\scriptstyle x_\pi, y_e$} (m-5-3)
		(m-7-5) edge  node[above] {$\scriptstyle x_d, y_e$} (m-5-3)
		(m-8-1) edge  node[left]  {$\scriptstyle \subset \{x_\pi, y_e\}$} (m-6-1)
		(m-8-1) edge (m-7-3)
		(m-9-3) edge  node[left] {$\scriptstyle x_d$} (m-7-3)
		(m-9-3) edge  node[above]  {$\scriptstyle \subset \{x_d\}$}  (m-8-1)
		% (m-8-1) edge  node[above]  {$\scriptstyle x_d$}  (m-7-3)
		;
		\path[->] 
		(m-9-3) edge  node[above] {$\scriptstyle x_\pi$} (m-7-5)
		(m-11-3) edge  node[left] {$\scriptstyle x_\delta$}  (m-9-3)
		(m-4-1) edge (m-2-2)
		(m-2-2) edge (m-1-3)
		;
		) ;
	\end{tikzpicture} 
	\caption{Some subgroups of the Monster $\mathbb{M}.$ }  
	\label{figure:subgroups:monster2}    
\end{figure}
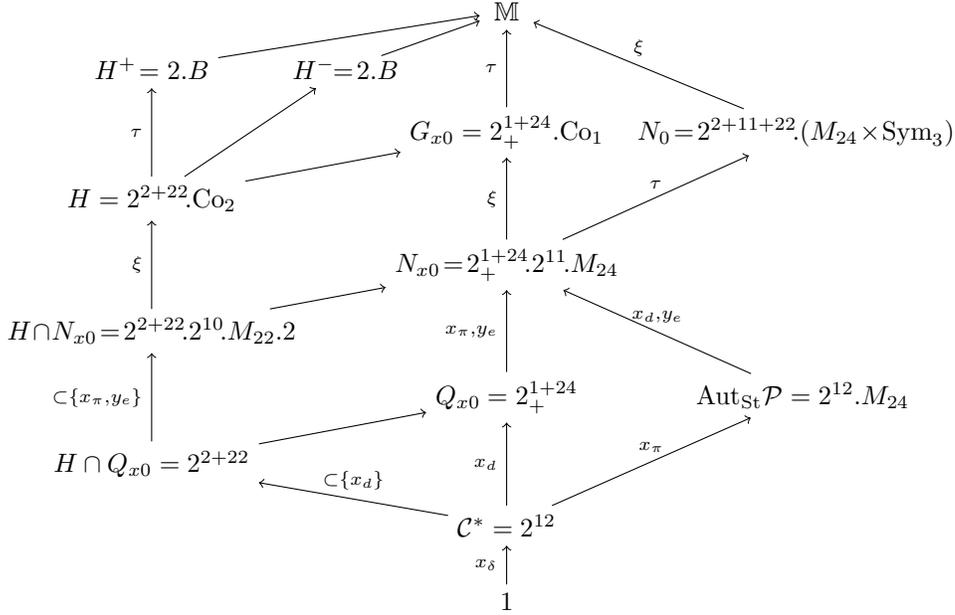

The subgroups of $\mathbb{M}$ relevant for our purposes (together with 
their structure) are shown
in Figure~\ref{figure:subgroups:monster2}. There arrows mean inclusion
of groups. If an arrow is labelled with one or more generators then all
generators of that shape must be added to the subgroup in order to 
obtain the group to which the arrow points. If an arrow is labelled 
with a set of generators preceded by a '$\subset$' symbol then some, but
not all generators of that shape must be added instead.

Note that $M_{22}$ is the Mathieu group acting on the set 
$\tilde{\Omega} \setminus \{2,3\}$ of size 22.

\subsection{Using axes of 2A involutions}
\label{sect:axes:2A}

If the centralizer $C_\mathbb{M}(t)$ of an element $t$ of $\setM$
fixes a unique one-dimensional subspace of the representation 
$196883_x$ of $\setM$ then a nonzero vector in that subspace is 
called an {\em axis} of $t$ in \cite{Conway:Construct:Monster}. 
Since $\rho = 196883_x +  1_x$, such an axis is unique 
in $\rho$ up to a scalar factor, and up to 
an additive multiple of the element $1_\rho$ in $1_x$.

Every 2A involution $t$
in $\mathbb{M}$  has an axis in $\rho$. Define vectors $v^+, v^-$ 
in in the subspace $300_x  \oplus 98280_x$ of $\rho$ by:
\[
v^+ = \eta_{2,2} + \eta_{3,3} -\eta_{2,3} - 2 X_\beta \, , \; 
v^- = \eta_{2,2} + \eta_{3,3} -\eta_{2,3} + 2 X_\beta \, ,
\]
with $\eta_{i,j}$ and $X_r$ as in Section~\ref{sect:rep:M}.
Then $v^+$ and $v^-$ are axes of   $x_{\beta}$ and  $x_{-\beta}$,
respectively, see \cite{Conway:Construct:Monster}. The axes 
$v^\pm$ have squared norm $\|v^\pm\| = 8$ and scalar product 
$\left< v^\pm , 1_\rho \right> = 2$. Also, $v^\pm *v^\pm$ is 
a positive scalar multiple of $v^\pm$, where the
operation '$*$'  $: \rho \times \rho \rightarrow \rho $ denotes
the Griess algebra defined in \cite{Conway:Construct:Monster}.
For any 2A involution  $t \in \mathbb{M}$ there is a unique axis 
satisfying these three conditions, and we write $\mbox{ax}(t)$ 
for that specific axis. In this paper an {\em axis} is always 
equal to $\mbox{ax}(t)$ for a 2A involution $t \in \setM$. 
For practical computations it suffices to know the co-ordinates
of an axis modulo 15, i.e. we may assume 
$\mbox{ax}(t) \in \rho_{15}$.

Let $v_1$ be as in Section~\ref{sect:Recogniz:Gx0}.
Given a triple $(v^+  g, v^-  g, v_1 g)$ of vectors in
$\rho_{15}$ for some unknown element $g$ of  $\mathbb{M}$ we want 
to find an element $h$ of  $\mathbb{M}$ (given as a word 
in $\Gamma^*$) with $g h = 1$. 

In Section \ref{sect:Orbits:Gx0} we construct a $h_1 \in \mathbb{M}$ 
that maps $v^+ \cdot g$ to $v^+$. Then $g h_1\in H^+$.

We call an axis in $\rho$ {\em feasible} if it is in the set 
$\{v^- h \mid h \in H^+\}$.
Since $g h_1 \in H^+$, the axis $v^- g h_1$ is feasible.
In Section~\ref{sect:Orbits:H} we construct a $h_2 \in H^+$ that 
maps the feasible axis $v^- \cdot g h_1$ to $v^-$. 
Note that $h_2$ fixes $v^+$. Therefore $g h_1 h_2$ 
fixes both, $v^+$ and $v^-$; hence 
$g h_1 h_2 \in H^+ \cap H^- \subset G_{x_0}$. Using $v_1 g$ and
the method in Section~\ref{sect:Recogniz:Gx0} we may find a 
$h_3 \in G_{x0}$ with  $g h_1 h_2 h_3 = 1$.

%%%%%%%%%%%%%%%%%%%%%%%%%%%%%%%%%%%%%%%%%%%%%%%%%%%%%%%%%%%%%%%%%%%%%%%%%%%%%%
\section{Reducing an axis of a 2A involution in the Monster}
\label{sect:Orbits:Gx0}
%%%%%%%%%%%%%%%%%%%%%%%%%%%%%%%%%%%%%%%%%%%%%%%%%%%%%%%%%%%%%%%%%%%%%%%%%%%%%%

In this section we show how to transform an arbitrary axis $v$
in $\rho$ to the standard axis $v^+$. The total number of axes 
in $\rho$ (or of 2A involutions in $\setM$) in is a bit less 
than $10^{20}$; so we cannot expect this task to be easy. 
The group $G_{x0}$ has 12 orbits on the set of 2A involutions, 
see \cite{Norton98}; and the co-ordinates of an axis of a 2A
involution encode quite a bit of geometric information about
that orbit. We will use this information
to find an $h \in \setM$ with $v h = v^+$.

For an implementation demonstrating that transformation of an axis we
refer to function {\tt reduce\textunderscore axis}
in \cite{mmgroup_doc},
Section {\em Demonstration code for the reduction algorithm}.

\subsection{Using the orbits of $G_{x0}$ on the axes for reducing an axis}
\label{ssect:Reducing:M}

For $v \in \rho$ let $M(v)$ be the real 
symmetric $24 \times 24$ matrix corresponding to the projection
of $v$ onto $300_x$, as in the Section~\ref{sect:Recogniz:Gx0}. 
An element $g$ of $G_{x0}$ acts as an orthogonal  transformation 
matrix $M_g$ on the Leech lattice, which is defined up to sign. 
Then $g$ maps $M(v)$ to $M_g^\top \cdot  M(v) \cdot M_g$. Hence
the eigenvalues of  $M(v)$ are actually properties of the orbit
$v G_{x0}$; so they can be used to identify the orbit 
$v G_{x0}$.

A transformation in $\mathbb{M}$ that maps an axis 
$v$ to $v^+$ has a representation as a word
\[
 g_1 \cdot \tau_1 \cdot g_2 \cdot \tau_2 \cdot
 \ldots \cdot  \tau_{n-1} \cdot  g_n  \; , 
\]
with $g_\nu \in G_{x0}$ and $\tau_\nu$ a power of $\tau$. 
Let $v$ be an axis.
Then  $\tau_\nu$ may change the orbit $vG_{x0}$, but $g_\nu$ 
does not. Let $l(v)$ be the minimum number of
occurrences of a power of $\tau$ in a word $h'$ of generators 
in $\Gamma^*$  such that $v h' = v^+$. 
Actually, $l(v)$ depends on the orbit $vG_{x0}$ only. 
So given an axis $v$ with $l(v) > 0$, we essentially have to
find a $h \in  G_{x0}$ such that $l(v \cdot h \tau^k) <l(v)$ 
for a $k = \pm 1$.

For an axis $v$, how does a power of $\tau$ change the 
orbit $vG_{x0}$? From (\ref{eqn:tau}) we conclude
\begin{equation}
   \label{eqn:orbit_Nx0}	
   N_{x0}  \tau^k    \subset \tau^k  N_{x0} \cup \tau^{-k} N_{x0}   
      \, , \quad \mbox{for} \; k = \pm 1 \; .
\end{equation}
Thus the set $\{v \tau^k G_{x0} \mid k = \pm1\}$ of 
orbits depends on the  orbit $v N_{x0}$ of $v$ only.
So given an axis $v$ we may first identify its orbit  $v G_{x0}$,
and then search for an orbit $v h N_{x0}$, $h \in G_{x0}$, such 
that  $\min\{ l(v h N_{x0}\tau^k) \, | \, k = \pm1 \} < l(v)$.

The group $N_{x0}$ is the stabilizer of the standard frame 
$\lambda_\Omega$ 
in $\Lambda / 2 \Lambda$ with respect to the action of $G_{x0}$.
Thus searching for a suitable $hN_{x0}$ in $G_{x0}$ amounts to 
searching for a type-4 vector $\lambda_r$ in $\Lambda / 2 \Lambda$ 
such that $\lambda_r h = \lambda_\Omega$. In
Appendix~\ref{app:map:type4} we will show how to
compute a representative of the coset  $hN_{x0}$ in  $G_{x0}$ from a
vector $\lambda_r \in \Lambda/2\Lambda$, such that we have
$\lambda_r h = \lambda_\Omega$. 
Note that both, $\lambda_\Omega$ and $\lambda_r$ are type-4
vectors in $\Lambda / 2 \Lambda$. So given an axis $v$ with
$l(v) > 0$, our task is to find a type-4 vector 
$\lambda_r \in \Lambda / 2 \Lambda$, such that for the coset  
$h N_{x0}$ with  $\lambda_r h = \lambda_\Omega$ we have 
$\min\{ l(v hN_{x0}\tau^k) \, | \,  k = \pm1 \} < l(v)$.
As we will see in the Section~\ref{ssect:Reducing:2A}, we 
may quickly identify the orbit $v G_{x0}$ of any  axis $v$. 
So finding an exponent $k = \pm1$ that minimizes $l(vh\tau^k)$ 
for a known $vh$ is easy.

For each axis $v$ with  $l(v) > 0$ we will specify a
reasonably small set $U_4(v)$ of type-4 vectors in $\Lambda/2\Lambda$
that can be used for decrementing the value $l(v)$ for an axis $v$
as described above. Here we have to treat the different orbits 
of $G_{x0}$ on the axes separately. Details are given in 
Section~\ref{ssect:mapping_axis}--\ref{ssect:Reducing:2A}.
Lemma~\ref{Lemma:orbit:2A} deals with axes $v$ satisfying $l(v)=0$.

\subsection{Enumeration of the orbits of $G_{x0}$ on the axes}
\label{ssect:Enumeration:Gx0}

Let  $v^+ = \mbox{ax}(x_\beta)$ be the axis of the 2A involution 
$x_\beta$ as defined in Section~\ref{sect:axes:2A}.
In this section we assume that a given vector $v$ in $\rho$ is
the axis $\mbox{ax}(t)$ of an (unknown) 2A involution~$t$. 

Norton \cite{Norton98} has enumerated the twelve orbits of $G_{x0}$ on
the 2A involutions in $\mathbb{M}$. For a 2A involution $t$ the 
orbit  of $t$ under the action of ${G_{x0}}$
is characterized by the of class $tx_{-1} $ in  $\mathbb{M}$, where 
$x_{-1}$ is the central involution in  $G_{x0}$.  The classes of 
$tx_{-1} $  (in ATLAS notation) are given in column~1 of
Table~\ref{table:orbit_Gx0_2A}. The powers of these classes are	
given in column~2. The structure of the centralizer of $tx_{-1} $  
in $G_{x0}$ is given in column~3, as in \cite{Norton98}.

We have computed representatives of the 12 orbits of $G_{x0}$ on
2A involutions with the software package \cite{mmgroup2020}. 
Here the essential information
required from  \cite{Norton98} is that there are no more than
12 such orbits. Then we may easily find such representatives by
computing images $v^+h$ for random elements $h$ of $\mathbb{M}$, 
and watermark them e.g. with the set of eigenvalues of the matrix 
$ M(\mbox{ax}(t)) = M(v^+ h)$, where $t = h^{-1}x_\beta h$. These 
eigenvalues are given in column~4 of Table~\ref{table:orbit_Gx0_2A}.
Column~5 contains the squared norm  $\|M(\mbox{ax}(t))\|$ of the 
matrix (modulo 15), which we define as the sum of the squares of 
the eigenvalues or, equivalently,  as the 
sum of the squares of the entries of  the symmetric matrix.

The package \cite{mmgroup2020} contains a function for computing the 
real character of the representation $196883_x$ of an element of
$\mathbb{M}$ that powers up to a 2B involution. Thus for any 2A 
involution $t$ we may check the class containing $tx_{-1}$ against
column~1 of Table~\ref{table:orbit_Gx0_2A}.

\begin{table}[h]
\centering	
\begin{tabular}{|c|c|c|c|c|}
		\hline
		$tx_{-1} $ & Powers & $C_{G_{x0}}(tx_{-1})$ & Eigenvalues 
		               & $\| M(\mbox{ax}(t)) \|$ \\
				   & of $tx_{-1} $  &  & of $256  M(\mbox{ax}(t))$  & 
				    $ \bmod 15 $ \\
		\hline
		2A   & & $2^{2+22}$.Co$_2$ &  $512^1, 0^{23}$ & 4 \\
		2B & & $2^{2+8+16}O_8^+(2)$ & $ 64^8, 0^{16}$ & 8 \\
		4A & 2B & $2^{1+22}.M_{23}$ &  $ 144^1, 16^{23}$ & 14  \\
		4B & 2A &  $(2^7 \times 2^{1+8}).S_6(2)$ & $72^1, 24^{16}, 8^7$ & 13 \\
		4C & 2B & $2^{1+14+5}.A_8$ & $32^8, 16^{16}$ & 3 \\
        6A & 2A & $2^2.U_6(2).2$ & $90^1, 26^1, 18^{22}$ & 4 \\
        6C & 2B & $2^{3+8}.(3 \times U_4(2)).2$ & $34^6, 18^{16}, 10^2$ & 5 \\
        6F & 2B  & $  2^{1+8}.A_9 $ & $24^{16}, 16^8$ &  14\\
        8B &  4A,2B  & $ 2.2^{10}.M_{11}  $ & $   48^1, 24^{12}, 16^{11} $ & 2 \\
		10A & 2A & 2.HS.2 & $60^1, 20^{22}, 12^{1} $ & 4 \\
		10B & 2B  & $ 2^{1+8}.(A_5 \times A_5).2 $ & 
		       $(24\pm 4 \sqrt{5} )^4, 20^{16}$  & 8 \\
        12C & 6A,4B,2A & $ 2 \times S_6(2) $ & $ 42^1, 26^7, 18^{16}$ &  10\\
		\hline 
\end{tabular}
\caption{Orbits of $G_{x0}$ on 2A involutions in $\mathbb{M}$}
\label{table:orbit_Gx0_2A}	
\end{table}

For labelling an orbit $\mbox{ax}(t) \, G_{x0}$  of $G_{x0}$ on the
set of the  axes we use the name of the corresponding class 
$t x_{-1}$ in column~1 of  Table~\ref{table:orbit_Gx0_2A};
and we put this name in single quotes. 
In Section~\ref{ssect:Reducing:2A} we will show how to
identify the orbit  $v G_{x0}$ of a given axis $v$.

Let $l(v)$ be as Section~\ref{ssect:Reducing:M}. In case 
$l(v)=0$ the  axis $v$ is in the orbit $v^+G_{x0}$ with name~'2A'.
From the definition of $v^+ = \mbox{ax}(x_\beta)$ in
Section~\ref{ssect:central:H} we see that matrix 
$M(v^+)$ has rank 1 and an eigenvalue 2, and that $M(v^+)$
is proportional to $\lambda_\beta^\top \otimes \lambda_\beta$. 
So we obtain:

\begin{Lemma}
	\label{Lemma:orbit:2A}
	The orbit '2A' of $G_{x0}$ on the axes consists of the elements
	$x_r$ of $Q_{x0}$ such that $\lambda_r = \lambda(x_r)$ is of type~2. 
	For any such  $x_r$ the kernel  
	$\ker (M(\mbox{ax}(x_r)) - 2 \cdot 1_\rho) $ is one-dimensional and
	spanned by the set $\Lambda(\lambda_r)$ of the shortest preimages 
	of $\lambda_r$ in the Leech lattice $\Lambda$.
\end{Lemma}

\subsection{Mapping an axis to a different orbit of $G_{x0}$}   
\label{ssect:mapping_axis}

For any  axis $v$ not in the orbit '2A' we will show how to 
compute a non-empty set $U_4 = U_4(v)$ of type-4 vectors
in $\Lambda / 2\Lambda$  with:
\begin{equation}
\label{eqn:U4}
\forall \lambda_r \in U_4(v)  \; \forall h  \in G_{x0} \, :
\lambda_r h  = \lambda_\Omega  \implies 
\min \{ l(vh \tau^{k})  \mid k = \pm 1 \} < l(v) \; .
\end{equation}
In Section \ref{ssect:Reducing:2A} we will construct a set
$U = U(v) \subset \Lambda / 2\Lambda$ such that the mapping 
$v \mapsto  U(v)$ visibly commutes with the action of $G_{x0}$; and 
we define  $U_4(v) $ to be the set of  type-4 vectors in $U(v)$.
So it suffices to check Equation~\ref{eqn:U4} for one representative
$v$ of each orbit of $G_{x0}$ on the  axes.  
By (\ref{eqn:orbit_Nx0}) it suffices to check each case 
of $v$ and each $\lambda_r \in U_4(v)$ for one $h  \in G_{x0}$ with  
$\lambda_r h  = \lambda_\Omega $.
This check can be done (and has been done) computationally.

By a sequence of such reduction steps we may reduce $l(v)$ to 0
as shown in  Figure~\ref{figure:Gx0:orbits:2A}.  
\begin{figure}[h]
	\centering
	\begin{tikzpicture}
	\matrix (m)
	[matrix of  nodes,row sep=2.0ex,column sep=2em,minimum width=2em]
	{  
		$l(v) \quad =$ 
		&    3  &   2  &   1  & 0   \\
		&  '6F' & '4C' \\
		& '10B' & '4B' & '2B'  \\
		& '12C' & '6A' &      & '2A'  \\
		& '10A' & '6C' & '4A' \\
		&       & '8B' \\
	};         	
	\path[->] 
	(m-2-2) edge (m-2-3) 
	(m-3-2) edge (m-2-3) 
	(m-2-3) edge (m-3-4) 
	(m-3-2) edge (m-3-3) 
	(m-3-3) edge (m-3-4) 
	(m-4-2) edge (m-3-3)
	(m-3-4) edge (m-4-5)
	(m-4-2) edge (m-4-3)
	(m-5-2) edge (m-4-3)
	(m-4-3) edge (m-5-4)
	(m-5-3) edge (m-5-4)
	(m-5-4) edge (m-4-5)
	(m-6-3) edge (m-5-4)
	;
	) ;
	\end{tikzpicture} 
	\caption{Reduction of an axis $v$ depending on its orbit $vG_{x0}$}  
	\label{figure:Gx0:orbits:2A}    
\end{figure}

In the remainder of this subsection we define certain subsets
of $\Lambda/ 2\Lambda$ depending on an axis $v$ that will be
used for constructing $U(v)$. 

We write $\rank_3(v,k)$ for the rank of the matrix
$M(v) - k \cdot 1_\rho$, considered as a matrix in
$\setF_3^{n \times n}$ by taking the entries of $M$ modulo~3.
This is well defined since $M(v)$ has entries in 
$\setZ[\frac{1}{2}]$. 
For an endomorphism $\phi$ of  $\setR^{24}$, with 
$\Lambda \subset \setR^{24}$,  we put
\[
\ker_{3,\Lambda} (\phi)  = 
(\phi^{-1}(3 \Lambda) \cap \Lambda) + 3 \Lambda \;. 
\]
This defines the kernel of $\phi$ with respect to the Leech 
lattice modulo $3$; and this kernel is a subspace of the space
$\Lambda / 3 \Lambda$ considered as a vector space over 
$\setF_3^{24}$. 

Matrix $M(v)$ acts naturally as an endomorphism of
$\setR^{24}$, with $\Lambda \subset \setR^{24}$.
For a 2A axis $v$ and an integer $k$ we define 
$\ker_3(v, k)$ by
\[ 
\ker_3(v, k) = 
  \{ w \in \ker_{3,\Lambda} (M(v) - k \cdot 1_\rho) 
  \mid 2 \leq \mbox{type(w)} \leq 4 \} + 2 \Lambda 
   \subset \Lambda / 2\Lambda \; . 
\]
The term $2 \Lambda$ in that formula implies 
$\ker_3(v, k) \subset \Lambda / 2 \Lambda$.
We need  $\ker_3(v,k)$ in case $\rank_3(v,k) = 23$ only. In this 
case we have $|\ker _3(v, k)| \leq 1$; and $\ker _3(v, k)$ can
easily be computed from matrix $M(v)$ with entries taken modulo~3.
For this computation we may use the information in
Table~\ref{table:orbits:Nx0}  about the 'shapes' of vectors 
in $\Lambda$ of type $\leq 4$.    
 
By construction, $\rank(v, 3)$ is a property of the orbit 
$v G_{x0}$; and the mapping $v \mapsto \ker_3(v, k)$ commutes with
the action of $G_{x0}$ on $300_x$ and $\Lambda$.

Apart from analysing matrix $M(v)$ we also consider the component
of the axis $v$ in the subspace $98280_x$ of $\rho$. $G_{x0}$ acts
monomially on  $98280_x$, without changing the number of 
co-ordinates with the same absolute value. For $v \in \rho_{15}$ 
we let $S_k(v)$ be the set of type-2 vectors $\lambda_r$ in 
$\Lambda/ 2 \Lambda$ such that the co-ordinate of $v$ with respect 
to the basis vector of  $98280_x$ corresponding to $\lambda_r$ has 
absolute value $k \pmod{15}$. 

If $S$ is a subset of $\Lambda / 2\Lambda$ then we write $\Span(S)$
for the linear subspace of $\Lambda / 2\Lambda$ generated by $S$,
and we write $S^\perp$ for the orthogonal complement of $S$ (with
respect to the natural scalar product) in $\Lambda/2\Lambda$.
Let $\rad(S) = \Span(S) \cap \Span(S)^\perp$ be the {\em radical}
of $S$; so the radical of a linear subspace $S$ of 
$\Lambda / 2\Lambda$ is the intersection of $S$ with its
orthogonal complement.

For any axis $v$ we construct $U(v)$ as a sum of certain sets
defined in this subsection, depending on the orbit $v G_{x0}$.
We write $U' + U''$ for the sum of two subsets  $U',  U''$
of $\Lambda / 2 \Lambda$.

For identifying the orbit $vG_{x0}$ of $v$ we sometimes evaluate
$M(v)$ at a type-2 vector $\lambda_r \in \Lambda/2\Lambda$. Define 
$M(v, \lambda_r) = w \cdot M(v) \cdot w^\top$, for a 
$w \in \Lambda(\lambda_r)$. This is well defined, since  
$\Lambda(\lambda_r)$ is unique up to sign if
$\lambda_r$ is of type~2.

\Skip{
 Here we will
analyse some of the small sets in that partition. 
One way to find suitable type-4 vectors is to
check the subspace of $\Lambda/2\Lambda$ generated by such a 
small set  for type-4 vectors. We may also 
analyse the radical of such a subspace, i.e. the intersection of 
that subspace with its orthogonal complement in  $\Lambda/2\Lambda$. 
We will see that it suffices to do our calculations in the 
representation $\rho_{15}$.
Computations that require divisions, such as the computation of
the kernel of a matrix, will be done in  $\rho_{3}$.
}

\subsection{Reducing an axis inside an orbit of $G_{x0}$}
\label{ssect:Reducing:2A}

In this subsection we assume that an axis $v$ is given.
For each of the 12 orbits of $G_{x0}$ on the  axes we will 
show how to check that axis $v$ is in that orbit. For each 
orbit different from  orbit  '2A' we will also show how to 
compute a set $U = U(v)$ such that the (non-empty) set $U_4(v)$
of type-4 vectors in $U(v)$ satisfies  (\ref{eqn:U4}).

\vspace{1ex}
\noindent {\bf Case '2A'}
\nopagebreak

We check $\|M(v)\|= 4 \pmod{15}$, $\rank_3(v, 2) = 23$; and we check
that  $\ker_3(v, 2)$ is a singleton $\{\lambda_r\}$ with 
$\lambda_r$ of type~2 and $M(v, \lambda_r) = 4 \pmod{15}$ .
We put  $U=\{\}$ in this case.

By Lemma~\ref{Lemma:orbit:2A} we have $v = \mbox{ax}(x_{\pm r})$ for an
$x_r$  with $\lambda(x_r) = \lambda_r$.
We can find a $h  \in G_{x0}$ that maps $\lambda_r$ to 
$\lambda_\beta$, as explained in Appendix~\ref{app:map:type2}.
Then $v h  = \mbox{ax}(x_{\pm \beta}) \in \{v^+, v^-\}$.
In case  $v h  = v^+$ we are
done. In case  $v h  = v^-$ we replace $v h $ by $v h  x_d$ for
a $d \in \mathcal{P}$ with $\langle d, \beta \rangle = 1$.

\vspace{1ex}
\noindent {\bf Case '2B'}
\nopagebreak

We check  $\|M(v)\|= 8 \pmod{15}$ and $\rank_3(v, 0) = 8$.
We have $|S_4(v)| = 120$. Put $U = \Span(S_4(v))$. Then $|U| = 256$; 
and $U$ contains 135 type-4 vectors. If $h  \in G_{x0}$ 
maps a type-4 vector in $U$ to $\lambda_\Omega$ then one of 
the axes  $v h  \tau^{\pm1}$ is in orbit '2A' of $G_{x0}$.

\vspace{1ex}
\noindent {\bf Case '4A'}
\nopagebreak

We check $\|M(v)\|= 14 \pmod{15}$ and  $\rank_3(v, 0) = 23$. 
Put $U = \ker_3(v, 0)$. Then $U$ is a singleton $\{\lambda_r\}$ with 
$\lambda_r$ of type~4. If $h  \in G_{x0}$ maps $\lambda_r$ 
to $\lambda_\Omega$ then one of the axes  $v h  \tau^{\pm1}$
is in orbit '2A' of $G_{x0}$.

\vspace{1ex}
\noindent {\bf Case '4B'}
\nopagebreak

We check $\|M(v)\|= 13 \pmod{15}$.
We have  $|S_1(v)| = 512$. Put $U = \rad(S_1(v))$. Then
$|U| = 128$; and $U$ contains 63 type-4 vectors.
If $h  \in G_{x0}$ maps a type-4 vector in $U$
to $\lambda_\Omega$ then one of the axes $v h  \tau^{\pm1}$ 
is in orbit '2B' of $G_{x0}$.

\vspace{1ex}
\noindent {\bf Case '4C'}
\nopagebreak

We check  $\|M(v)\|= 3 \pmod{15}$. 
We have  $|S_1(v)| = 16$.  Put $U = \rad(S_1(v))$. Then
$|U| = 32$; and $U$ contains 15 type-4 vectors. 
If $h  \in G_{x0}$ maps a type-4 vector in $U$
to $\lambda_\Omega$ then one of the axes $v h  \tau^{\pm1}$ 
is in orbit '2B' of $G_{x0}$.

\vspace{1ex}
\noindent {\bf Case '6A'}
\nopagebreak

We check $\|M(v)\|= 4 \pmod{15}$, $\rank_3(v, 2) = 23$; and we check
that  $\ker_3(v, 2)$ is a singleton $\{\lambda_r\}$ with 
$\lambda_r$ of type~2 and $M(v, \lambda_r) = 7 \pmod{15}$ .
Put $U = \{\lambda_r\} + S_5(v)$. Then $|U| = 891$; and all vectors
in  $U$ are of type~4.
If $h  \in G_{x0}$ maps one vector in $U$ 
to $\lambda_\Omega$ then one of the axes  $v h  \tau^{\pm1}$
is in orbit '4A' of $G_{x0}$.

\vspace{1ex}
\noindent {\bf Case '6C'}
\nopagebreak

We check $\|M(v)\|= 5 \pmod{15}$. 
We have  $|S_3(v)| =36$; Put $U = \Span(S_3(v))$ Then
$|U| = 64$; and  $U$ contains 27 type-4 vectors. 
If $h  \in G_{x0}$ maps one type-4 vector in $U$ 
to $\lambda_\Omega$ then one of the axes  $v h  \tau^{\pm1}$
is in orbit '4A' of $G_{x0}$.

\vspace{1ex}
\noindent {\bf Case '6F'}
\nopagebreak

We check $\|M(v)\|= 14 \pmod{15}$ and $\rank_3(v, 0) = 8$.
We have  $|S_7(v)| =144$.  Put $U = \rad(S_7(v))$. Then
$|U| = 256$; and $U$ contains 135 type-4 vectors.
If $h  \in G_{x0}$ maps one type-4 vector in $U$ 
to $\lambda_\Omega$ then one of the axes  $v h  \tau^{\pm1}$
is in orbit '4C' of $G_{x0}$.

\vspace{1ex}
\noindent {\bf Case '8B'}
\nopagebreak

We check $\|M(v)\|= 2 \pmod{15}$. 
We have  $|S_1(v)| =24$. Put $U = S_1(v) + S_1(v)$.
Then $|U| = 134$; and $U$ contains one type-4 vector
$\lambda_r$. If $h  \in G_{x0}$ 
maps $\lambda_r$ to $\lambda_\Omega$ then one of the axes  
$v h  \tau^{\pm1}$ is in orbit '4A' of $G_{x0}$.
Note that $\lambda_r \in \{\lambda_s\} + S_1(v)$ holds for all
$\lambda_s \in S_1(v)$.

\vspace{1ex}
\noindent {\bf Case '10A'}
\nopagebreak

We check $\|M(v)\|= 4 \pmod{15}$ and $\rank_3(v,2) = 2$.   
Then $|S_3(v)|=1$ and $|S_1(v)| = 100$.  Put  
$U = S_3(v) + S_1(v)$.
Then all the 100 vectors in $U$ are of type 4.
If $h  \in G_{x0}$ maps one  vector in $U$ 
to $\lambda_\Omega$ then one of the axes  $v h  \tau^{\pm1}$
is in orbit '6A' of $G_{x0}$.

\vspace{1ex}
\noindent {\bf Case '10B'}
\nopagebreak

We check $\|M(v)\|= 8 \pmod{15}$ and $\rank_3(v,0) = 24$.
Then  $|S_4(v)| = 680$. Put $U = \rad(S_4(v))$.
Then $|U| = 256$; and $U$ contains 135 type-4 vectors. 
If $h  \in G_{x0}$ maps one type-4 vector in $U$  to
$\lambda_\Omega$ then one of the axes  $v h  \tau^{\pm1}$ 
is in the orbit '4B' or '4C'  of $G_{x0}$.
An axis can be mapped to orbit '4B' (or '4C') for
60 (or 75) of the 135 type-4 vectors, respectively.

\vspace{1ex}
\noindent {\bf Case '12C'}
\nopagebreak

We check $\|M(v)\|= 10 \pmod{15}$. 
We have  $|S_7(v)| = 56$. Put $U = \rad(S_7(v))$.
Then $|U| = 256$; and $U$ contains 135 type-4 vectors. 
If $h  \in G_{x0}$ maps one type-4 vector in $U$  to
$\lambda_\Omega$ then one of the axes  $v h  \tau^{\pm1}$ 
is in the orbit '4B' or '6A'  of $G_{x0}$.
An axis can be mapped to orbit '4B' (or '6A') for
63 (or 72) of the 135 type-4 vectors, respectively. 

\vspace{1ex}
\noindent {\bf Remark}
\nopagebreak

It has not escaped our attention that the sets $U$ of  vectors
in $\Lambda / 2\Lambda$ discussed in the different cases for 
$v = \mbox{ax}(t)$ provide a wealth of geometric information 
about the centralizer of the 2A involution $t \in \setM$  in 
$G_{x0}$. According to Table~\ref{table:orbit_Gx0_2A}, in Case~'10A'
the Higman-Sims group HS is involved in that centralizer.
The 100 short vectors in $\Lambda / 2\Lambda$  mentioned in that 
case correspond to the 100 short vectors in $\Lambda$ discussed
in the analysis of the  Higman-Sims group in \cite{Conway-SPLG}, 
Ch.~10.3.5. In Case '10B' the sum of the two 4-dimensional 
eigenspaces of matrix $M(\mbox{ax}(t))$ contains a sublattice
$\sqrt{2} E_8$ of $\Lambda$; and this can be interpreted in terms 
of the {\em icosians}, see \cite{Conway-SPLG}, Ch.~8.2.1.
Most of the other cases have simpler geometric interpretations. 
As the  objective of this paper is computational, 
we omit further details.

%%%%%%%%%%%%%%%%%%%%%%%%%%%%%%%%%%%%%%%%%%%%%%%%%%%%%%%%%%%%%%%%%%%%%%%%%%%%%%
\section{Reducing a feasible axis in the group $H^+$}
\label{sect:Orbits:H}
%%%%%%%%%%%%%%%%%%%%%%%%%%%%%%%%%%%%%%%%%%%%%%%%%%%%%%%%%%%%%%%%%%%%%%%%%%%%%%

In this section we will define the term {\em feasible axis}; and we
will show how to transform an arbitrary feasible axis $v$
in $\rho$ to the standard feasible axis $v^-$ using a computation
in the group $H^+$. 
For an implementation demonstrating that transformation of an 
feasible axis we refer to function
{\tt reduce\textunderscore feasible\textunderscore axis}
in \cite{mmgroup_doc},
Section {\em Demonstration code for the reduction algorithm}.

\subsection{Enumeration of the orbits of $H$ on feasible axes}
\label{ssect:Enumeration:H}

Let $v^+, v^-$ be the axes of the 2A involutions $x_\beta, x_{-\beta}$, 
and let $H^+$ be the centralizer of $v^+$ (or of $x_\beta$) 
as in Section ~\ref{sect:strategy}.
We call an axis in $\rho$ {\em feasible} if it is in the set 
$\{v^- h \mid h \in H^+\}$, as in Section~\ref{ssect:central:H}.
Since $H^+$ fixes $v^+$, and  $v^+$ is orthogonal to  $v^-$
with respect to the scalar product given by the natural norm in 
$\rho$, any feasible axis is also orthogonal  to $v^+$. For two
2A involutions  $t_1, t_2 \in \mathbb{M}$
their axes $\mbox{ax}(t_1)$ and  $\mbox{ax}(t_2)$ are 
orthogonal if and only if the product $t_1\cdot t_2$ is 
a 2B~involution in $\mathbb{M}$; and  $\mathbb{M}$  is transitive
on the pairs of orthogonal axes; see \cite{Conway:Construct:Monster}
or \cite{citeulike:Monster:Majorana}. 

Given a feasible axis $v$, we construct a $g \in H^+$ with $v g$ = $v^-$
in this section. Therefore we have to compute in the group $H^+$. 
The group  $H^+ \cap N_{x0}$ is generated by
\[
   \Gamma_f^0 = 
   \left\{x_\delta, x_d, y_d, x_\pi \in \Gamma \mid
   \delta \in \mathcal{C}^*, d \in \mathcal{P}, 
   \langle d, \beta \rangle = 0, \pi \in \AutStP,  \beta^\pi = \beta
   \right\} \; ;
\]
so computation in $H^+ \cap N_{x0}$ is easy. The group 
$H = H^+ \cap G_{x0}$ is generated by $\Gamma_f^0$ and $\xi$; 
the group $H^+$ is generated by $H$ and  $\tau$. Put
$\Gamma_f = \Gamma_f^0 \cup \{\xi^{\pm1}, \tau^{\pm1}\}$,
so that $\Gamma_f$ generates $H^+$. Let $\Gamma_f^*$ be the
set of words in $\Gamma_f$.

Our strategy for finding a  $g \in H^+$ that maps a feasible axis
to $v^-$ is similar to the strategy in Section~\ref{sect:Orbits:Gx0}
for finding an element of $\setM$ that maps an arbitrary axis to the
standard axis $v^+$. Here we will replace the r\^oles of the groups 
$\mathbb{M}$, $G_{x0}$, and $N_{x0}$ by their intersections with the 
group $H^+$, as indicated in Table~\ref{reduce:analog}.

\begin{table}[!h]
	\centering	
	\begin{tabular}{l | c|c}
		\hline
		 Reduction operation &
		Reduce axis to $v^+$  &  Reduce feasible axis to $v^-$  \\
		\hline
		  Operating  in group &
		    $\mathbb{M}$   &  $ H^+  $ \\
		  Set of generators &
		  $\Gamma$ & $\Gamma_f$ \\  
		  Using orbits of group &
		   $G_{x0} $       &  $  H $ \\
		  Number of orbits on axes &  12, see\cite{Norton98} 
		    & 10, see \cite{mueller_2008} \\
		  Reduction path &  see Figure~\ref{figure:Gx0:orbits:2A}  
		        &   see Figure~\ref{figure:H:orbits:2A}   \\  
		  Vectors $\lambda_r  \! \in \! \Lambda/ 2 \Lambda$ used   &
           $\mbox{type}(\lambda_r) \! =  \! 4$  & 
           $\mbox{type}(\lambda_r) \! =  \!  2  , 
           \mbox{type}(\lambda_r \! + \! \lambda_\beta) \! =  \!  4$  \\
           $\lambda_r$ is mapped to image &  $\lambda_{\Omega}$ &
                $\lambda_{\Omega} + \lambda_\beta$    \\
           Group fixing that image  &
           $N_{x0} $       &  $   H \cap N_{x0} $ \\                
		\hline
	\end{tabular}
	\caption{Analogies between reducing axes to  $v^+$ and feasible axes to $v^-$}
	\label{reduce:analog}	
\end{table}

We present an element $g$ of $H^+$ as a word 
\[
g_1 \cdot \tau_1 \cdot g_2 \cdot \tau_2 \cdot
\ldots \cdot  \tau_{n-1} \cdot  g_n  \; , 
\]
with $g_\nu \in H $ and $\tau_\nu$ a power of $\tau$. 
Let $v$ be a feasible axis.
Then  $\tau_\nu$ may change the orbit $vH$, but $g_\nu$ does not.
Let $l(v)$ be the minimum number of occurrences of a power of 
$\tau$ in a word $h'$ in $\Gamma_f^*$  such that $v h' = v^-$. 
Actually, $l(v)$ depends on the orbit $vH$ only. So given a 
feasible  axis $v$ with $l(v) > 0$ we have to to find a $h \in H$
 such that $l(v \cdot h  \tau^k) <l(v)$ for a $k = \pm 1$.

Similar to equation \ref{eqn:orbit_Nx0}	we obtain
\begin{equation}
	\label{eqn:orbit_Nx0_cap_H}	
	(H \cap N_{x0})  \tau^k    \subset \tau^k 
	(H \cap N_{x0}) \cup \tau^{-k} (H  \cap N_{x0})   
	\, , \quad \mbox{for} \; k = \pm 1 \; .
\end{equation}
Thus the set $\{v \tau^k H \mid k = \pm1\}$ of 
orbits depends on the  orbit $v (H \cap N_{x0})$ of $v$ only.
So given an axis $v$ we may first identify its orbit  $v H$,
and then search for an orbit $v h  (H \cap N_{x0})$, 
$h  \in H$, such that  
$\min\{ l(v h  (H \cap N_{x0})\tau^k) \, | \, k = \pm1 \} < l(v)$.

We call a vector $\lambda_r$ in $\Lambda/ 2\Lambda$ {\em feasible} if
$\mbox{type}(\lambda_r) = 2$ and 
$\mbox{type}(\lambda_r + \lambda_\beta) = 4$. This is equivalent to
claiming that any two shortest preimages of $\lambda_r$ and  
$\lambda_\beta$ in the Leech lattice $\Lambda$ are of type~2 
and perpendicular. The group $H$ fixes $\beta$;  and it is transitive 
on the set of feasible vectors in $\Lambda / 2 \Lambda$. The group 
$H \cap N_{x0}$ fixes $\lambda_\beta$ and $\lambda_\Omega$, and 
hence also the feasible type-2 vector $\lambda_\Omega + \lambda_\beta$. 
From \cite{Atlas} and the structures of $H$ and of $H \cap N_{x0}$
we see that  $H \cap N_{x0}$ is a maximal subgroup of $H$. So
$H \cap N_{x0}$ is the stabilizer of the feasible vector  
$\lambda_\Omega + \lambda_\beta$ in $H$.

Thus searching for a suitable coset $h (H \cap N_{x0})$ in $H$ amounts 
to searching for a feasible vector $\lambda_r$ in $\Lambda / 2\Lambda$ 
such that $\lambda_r h  = \lambda_\Omega  + \lambda_\beta$. In
Appendix~\ref{app:map:feasible} we will show how to
compute a representative of the coset  $h (H \cap N_{x0})$  in $H$
from a feasible vector $\lambda_r$ satisfying  
$\lambda_r h  = \lambda_\Omega + \lambda_\beta$.

M\"uller \cite{mueller_2008} has enumerated the ten orbits of
$2^{1+22}.\mbox{Co}_2$ on the representation of the 
Baby Monster $B$ acting on the cosets of  $2^{1+22}.\mbox{Co}_2$ 
as a permutation group. 
Thus the group $H^+$ of structure $2.B$ has ten orbits of
$H$ (of structure $2^{2+22}.\mbox{Co}_2$) acting on the cosets 
of $H$. The centralizer of the feasible  axis $v^-$ 
in $H^+$ is $H^+ \cap H^- = H$. Since $H^+$ is transitive on the set of
feasible axes,  we conclude that there are ten orbits of $H$ on the
feasible axes.

We have computed representatives of the 10 orbits of $H$ on the
feasible axes with the software package \cite{mmgroup2020}. 
Here the essential information
required from  \cite{mueller_2008} is that there are no more than
10 such orbits. Then we may easily find such representatives by
computing images $v^-h$ for random elements $h$ of $H^+$.
Here we watermark the feasible axes in the same way as the
axes in Section~\ref{sect:Orbits:Gx0}; and we augment our 
watermark of an axis $v$ by the value $M(v, \lambda_\beta)$,
as defined in Section~\ref{ssect:mapping_axis}. 
$M(v, \lambda_\beta)$ is invariant under the action of $H$.

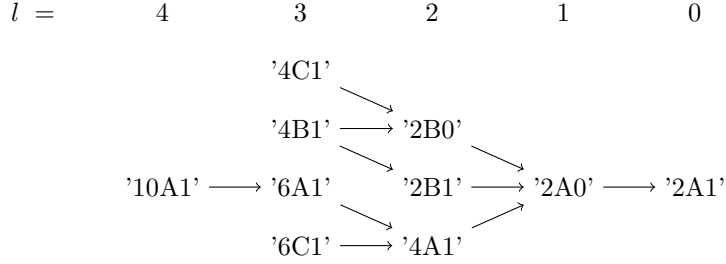
\begin{figure}[!h]
	\centering
	\begin{tikzpicture}
		\matrix (m) 
		[matrix of  nodes,row sep=2.0ex,column sep=2em,minimum width=2em]
		{   
			$l \; = $  &    4   &  3   &  2  &   1 & 0   \\
			           &        & '4C1'  \\
		          	   &        & '4B1'  &  '2B0'  \\
			           & '10A1' & '6A1'  &  '2B1'  & '2A0' & '2A1'  \\
		               &        & '6C1'  &  '4A1' \\
		};         	
		\path[->] 
		(m-2-3) edge (m-3-4) 
		(m-3-3) edge (m-3-4) 
		(m-3-3) edge (m-4-4)
		(m-3-4) edge (m-4-5)
	 	(m-4-4) edge (m-4-5)
	 	(m-4-5) edge (m-4-6)
		(m-4-2) edge (m-4-3)
		(m-4-3) edge (m-5-4)
		(m-5-3) edge (m-5-4)
		(m-5-4) edge (m-4-5)
		;
		) ;
	\end{tikzpicture} 
	\caption{Reduction of orbits of $H$ on feasible  axes}  
	\label{figure:H:orbits:2A}    
\end{figure}

For labelling  an orbit of $H$ on the set of the feasible  axes
we use the name of the orbit of  $G_{x0}$ on the set of all axes,
which contains that orbit of $H$. These names are defined as in  
Section~\ref{ssect:Enumeration:Gx0}. We append a '0' digit 
to that name if $M(v, \lambda_\beta) = 0$; otherwise we append 
a '1' digit. This labelling is sufficient to distinguish
between the 10 orbits.
Thus for identifying the orbit of $H$ containing the feasible axis
$v$ we first have to compute the orbit $v G_{x0}$ on the axes as 
described Section~\ref{ssect:Reducing:2A}. If there are several
orbits  $v H$ on the feasible  axes contained in the same
orbit $v G_{x0}$  then we also have to compute 
$M(v, \lambda_\beta)$. 
It suffices to compute $M(v, \lambda_\beta)$ modulo 15.
We obtain the 10 orbits of  $H$ on the feasible axes shown in
Figure~\ref{figure:H:orbits:2A}.   

\subsection{Reducing a feasible axis}
\label{ssect:Reducing:feasible}

In this subsection we assume that a feasible  axis $v$ is 
given. Let $U(v)$ be as in Section~\ref{ssect:mapping_axis}.
We define the set $U_f(v)$  of feasible type-2 vectors
in $\Lambda / 2 \Lambda$  by:
\begin{equation}
\label{eqn:Uf}
U_f(v) = \{\lambda_r + \lambda_\beta \mid \lambda_r \in U(v), \;
\mbox{type}(\lambda_r) = 4, \;
\mbox{type}(\lambda_r + \lambda_\beta) = 2\} \; , 
\end{equation}
If $v$ is in one of the orbits of $H$ on the feasible axes shown in 
Figure~\ref{figure:H:orbits:2A}	(except for orbits '2A0' and '2A1'),
then $U_f(v)$ is non-empty, and  we can show computationally:
\[
\forall \lambda_r \in U_f(v)  \; \forall h  \in H \, :
\lambda_r h  = \lambda_\Omega + \lambda_\beta  \implies 
\min \{ l(vh \tau^{k})  \mid k = \pm 1 \} < l(v) \; 
\]
Therefore we use a similar method as described in 
Section~\ref{ssect:mapping_axis}.

We obtain the following cases of orbits of $H$ on feasible axes.

\vspace{1ex}
\noindent {\bf Case '2A1'}
\nopagebreak

Since $v^- = \mbox{ax}(x_{-\beta})$  is in this orbit of $H$, and $H$
fixes $\lambda_\beta$, the axis $v^-$ is the only feasible axis in 
this orbit, so that we are done.

\vspace{1ex}
\noindent {\bf Case '2A0'}
\nopagebreak

By Lemma~\ref{Lemma:orbit:2A} there is an $x_r$ with 
$v = \mbox{ax}(x_r)$, and $x_r$ can be effectively computed from
$v$. Since $v$ is feasible and $M(v, x_\beta) = 0$,  the vectors
in $\Lambda(\lambda_r)$ and  $\Lambda(\lambda_\beta)$ 
in the Leech lattice are perpendicular, and hence 
$\mbox{type}(\lambda_\Omega + \lambda_r) =4$, i.e. $\lambda_r$ is
feasible. 

So we can find a $h  \in H$ that maps $\lambda_r$ to 
$\lambda_\Omega+\lambda_{\beta}$ without changing  $\lambda_\beta$, 
as explained in Appendix~\ref{app:map:feasible}.
Then $vh  = \mbox{ax}(x_{\pm \Omega}  x_\beta)$.
Note that the
element $\tau$ of $H^+$ cyclically permutes the axes
$\mbox{ax}(x_{-\beta}) $, $\mbox{ax}(x_{-\Omega} x_\beta) $,
and  $\mbox{ax}(x_{\Omega} x_\beta) $. Thus there is a $k = \pm 1$
with  $v h  \tau^k = \mbox{ax}(x_{-\beta}) = v^- $; i.e.
$v h  \tau^k$  is the unique axis in the orbit '2A1'.

\Skip{
\vspace{1ex}
\noindent {\bf Cases '2B0' and '2B1'}
\nopagebreak

The set $U_f(v)$ contains 15 feasible type-2  vectors in case '2B0' 
and 63 feasible  type-2 vectors in case '2B1'. 
For any $\lambda_s \in U_f(v)$ we can find a $h  \in H$ that maps 
$\lambda_s$  to  $\lambda_\Omega+\lambda_{\beta}$, as explained 
in Appendix~\ref{app:map:feasible}.
Then one of the axes $v h  \tau^{\pm1}$ is in the orbit '2A0'
of $H$.
}

\vspace{1ex}
\noindent {\bf The remaining cases
  %	'2B0', '2B1','4A1', '4B1', '4C1', '6A1', '6C1', and '10A1'
  of orbits of $H$ on feasible axes
}
\nopagebreak

The sizes of the sets $U_f(v)$ in these cases are given 
in Table~\ref{table:map:feasible}.
For any $\lambda_s \in U_f(v)$ we can find a $h  \in H$ that maps 
$\lambda_s$  to  $\lambda_\Omega+\lambda_{\beta}$, as explained 
in Appendix~\ref{app:map:feasible}.
Then one of the axes  $v h  \tau^{\pm1}$ is in an orbit of $H$ 
as shown in  Table~\ref{table:map:feasible}.

In case '4B1' the set $U_f(v)$ has size 31; for 30 elements of 
$U_f(v)$ one of the axes $vh  \tau^{\pm1}$ is in the orbit  '2B0';
for the remaining element of $U_f(v)$ one of these axes is 
in the orbit  '2B1'.

\begin{table}[!h]
	\centering	
	\begin{tabular}{| l |  c | c | c | c | c | c | c | c | c |   }
		\hline
		Orbit of $H$ & 2B0 & 2B1 &  4A1 &
		\multicolumn{2}{c|}{4B1}  &
		4C1 & 6A1 & 6C1 &  10A1 \\
		\hline
		is mapped to orbit &
	    2A0 & 2A0 & 2A0 & 2B0 & 2B1 & 2B0 & 4A1 & 4A1 & 6A1 \\
	    \hline
		Size of $U_f$ &
		15 & 63 & 1  & 30 & 1 & 15 & 891 & 27 & 100 \\
		\hline
	\end{tabular}
	\caption{Reduction of feasible axes}
	\label{table:map:feasible}	
\end{table}

% \newpage

%%%%%%%%%%%%%%%%%%%%%%%%%%%%%%%%%%%%%%%%%%%%%%%%%%%%%%%%%%%%%%%%%%%%%%%%%%%%%%
\section{Conclusion}
\label{sect:Conlsusion}
%%%%%%%%%%%%%%%%%%%%%%%%%%%%%%%%%%%%%%%%%%%%%%%%%%%%%%%%%%%%%%%%%%%%%%%%%%%%%%

We have shown how to reconstruct an element $g$ of the Monster
$\mathbb{M}$ as a word in the generators of $\mathbb{M}$  from 
the images of three fixed vectors in the representation $\rho$ 
under the action of $g$. It suffices if the co-ordinates of 
these three images are  known modulo~15. This leads to an extremely 
fast word shortening algorithm.
In \cite{mmgroup2020} we have implemented  the group operation
of the Monster $\mathbb{M}$ based on this word shortening algorithm. 

By construction, the new algorithm computes a unique representation 
of an element of $\mathbb{M}$. This representation can be compressed 
to a bit string of less than 256 bit length. So we may quickly find 
an element of $\mathbb{M}$ in an array of millions of
elements of $\mathbb{M}$.

The package \cite{mmgroup2020} can also perform the following
computations efficiently. 
Given two arbitrary conjugate involutions $t_1, t_2 \in \setM$,
we can compute a $g \in \setM$ with $t_1^g = t_2$.  
In the smallest faithful real representation $196883_x$ of 
$\mathbb{M}$ we can compute the character of an element 
of $\mathbb{M}$ centralizing a {\em known} 2B involution.
Therefore we also use a strategy called 'changing post'
in \cite{Wilson13}.

The ATLAS \cite{Atlas} describes a well-known homomorphism from 
a certain Coxeter group $Y_{555}$ to a group of structure 
$(\setM \times \setM).2$, which is called the {\em Bimonster}. 
For background, see \cite{ivanov_1999}, Section~8. 
The software package \cite{mmgroup2020} also contains the first 
efficient implementation of this homomorphism, based on the 
ideas in \cite{Far12, Nor02}. During the construction of that
homomorphism in \cite{mmgroup2020} we have verified the 
presentation of the Monster given in \cite{Nor02}.

The new algorithm has been used to solve the long-standing 
problem of finding all maximal subgroups of the Monster, 
see \cite{DLP2023}.

In principle, the new algorithm can  also compute unique 
representatives of the cosets in a chain of subgroups 
\[
\ldots \, \subset  \, 2^{2+22} \!.2^{10}\! .M_{22}.2
\, \subset \, 2^{2+22}\!.\mbox{Co}_2
\, \subset \, 2.B  \, \subset \, \mathbb{M} .
\]
So using the methods in
\cite{MNW_2006}, we may be able to enumerate rather big orbits of 
some very large subgroups of the Monster $\mathbb{M}$.

\fatline{Acknowledgements}: 
I would like to thank Heiko Dietrich for many helpful comments and 
suggestions on the paper. I would also like to thank Max Horn 
and Dima Pasechnik for valuable help on the continuous integration
of the software project.

 \newpage

%%%%%%%%%%%%%%%%%%%%%%%%%%%%%%%%%%%%%%%%%%%%%%%%%%%%%%%%%%%%%%%%%%%%%%%%%%%%%%
\section*{Notation}
\label{Notation}

\begin{tabular}{cll}
 Symbol & Description 
 \hphantom{01234567890123456789012345678901234567890123456789}
 & \hspace{-2.25em} Section \vspace{0.5ex} \\
$\AutStP$ &
The group of standard automorphisms of the Parker
loop $\mathcal{P}$.   
& \ref{ssect:Parker}  \\
 axis &
Vector in representation $\rho$, corresponding to
a 2A involution in $\setM$. 
& \ref{sect:axes:2A}  \\      
$\mbox{ax}(t) $ &
The  axis of a 2A involution $t \in \setM$
in the representation $\rho$.
& \ref{sect:axes:2A}  \\ 
$B$ &
The Baby Monster group. $\setM$ has a subgroup
of structure $2.B$.
&  \ref{ssect:central:H}\\
$\beta$ &
A fixed element of $\mathcal{C^*}$ used for constructing
$x_\beta, x_{-\beta}$, and $\lambda_\beta$.
& \ref{ssect:central:H} \\
$\mathcal{C}$, $\mathcal{C}^*$ &
$\mathcal{C}$ is the 12-dimensional Golay code in $\Field_2^{24}$;
$\mathcal{C}^*$ is its cocode  $\Field_2^{24}/\mathcal{C}$.
& \ref{ssect:Golay}  \\
$\mbox{Co}_0, \mbox{Co}_1$ &
Automorphism groups of $\Lambda$ and of $\Lambda / 2 \Lambda$.
$\mbox{Co}_0$ has structure $2.\mbox{Co}_1$.
&  \ref{ssect:Leech:Qx0}  \\  
$\mbox{Co}_2$ &
The subgroup of $\mbox{Co}_1$ fixing a shortest vector
in  $\Lambda / 2 \Lambda$.
&  \ref{sect:axes:2A}  \\  
$d,e,f$ &
    Elements of the Parker loop $\mathcal{P}$ or of
    the Golay code $\mathcal{C}$.
    &  \ref{ssect:Parker}  \\
$\delta, \epsilon$ &
Elements of the Golay cocode  $\mathcal{C}^*$.
& \ref{ssect:Parker}  \\
$\eta_{i} $ &
    $i$-th unit vector of the space $\setR^{24}$ containing $\Lambda$,
    for $i \in \tilde{\Omega}$.
    & \ref{ssect:Leech:Qx0} \\
$\eta_{i,j} $ &
Matrix in $300_x$ with entry 1 in row $i$, column $j$, 
and in row $j$, col. $i$. 
& \ref{sect:rep:M} \\
feasible & 
Property of an axis or a vector in $\Lambda / 2 \Lambda$,
related to involutions $x_{\pm \beta}$.
& \ref{ssect:Enumeration:H} \\
$G_{x0} $ &
Maximal subgroup of $\setM$ of structure $2^{1+24}_+\!.\mbox{Co}_1$
with centre $\{1, x_{-1}\}$. 
& \ref{ssect:maximal:G0}  \\ 
$\Gamma, \Gamma^*$ &   
$\Gamma$ is a  set of elements generating $\setM$. 
$\Gamma^*$ is the set of words in $\Gamma$.
& \ref{sect:rep:M}   \\
$H$ &
Group of structure $2^{2+22}\!.\mbox{Co}_2$, with
$H = H^+ \cap H^- = G_{x0} \cap H^+$.
& \ref{ssect:central:H}   \\
$H^+, H^-$ &
Centralizers of the involutions $x_\beta, x_{-\beta}$ in $\setM$,
both of structure $2.B$.
    & \ref{ssect:central:H} \\
$i,j$ &
Elements of $\tilde{\Omega}$,
also considered  as elements of  $\mathcal{C}^*$  of weight 1.
&  \ref{ssect:Parker}  \\       
$\Lambda$, $\Lambda/ 2\Lambda$ &
The 24-dimensional Leech lattice, and the Leech lattice mod 2.
& \ref{ssect:Leech:Qx0} \\
$\Lambda(\lambda_r) $ &
For  $\lambda_r \in\Lambda/ 2\Lambda $ this is the
set of shortest preimages of $\lambda_r$  in $\Lambda$.  
& \ref{ssect:Leech:Qx0} \\
$\lambda$ &
Homomorphism from $Q_{x0}$ onto $\Lambda/ 2\Lambda$ with
kernel $\{1, x_{-1}\}$.
& \ref{ssect:Leech:Qx0} \\
$ \lambda_{\beta}$ &
Fixed vector in $\Lambda/ 2\Lambda$ of type 2, equal to
$\lambda(x_{\beta})$. 
& \ref{ssect:central:H}   \\
$  \lambda_{\Omega}$, &
Fixed vector in $\Lambda/ 2\Lambda$ of type~4; 
$\mkern-1.5mu \Lambda \mkern-1mu   (\!\lambda_{\Omega}\!)$ 
is the co-ordinate frame of $\Lambda$.
& \ref{ssect:Leech:Qx0} \\
$\setM$ &
    The Monster group, i.e. the largest sporadic simple group.   
    &  \ref{sect:Introduction}  \\ 
$M_{24}$ &
Mathieu group,  acts on $\tilde{\Omega}$ as
the automorphism group of $\mathcal{C}$.
& \ref{ssect:Golay} \\  
$M(v)$ &
Projection of a vector $v \in \rho$ onto the subspace
$300_x$ of  $\rho$; also \\
& considered as a symmetric $24 \times 24$ matrix,
with  $\Lambda \subset \setR^{24}$.   
&\ref{sect:Recogniz:Gx0} \\  
$M(v, \lambda_r)$ &
Equal to $u M(v) u^\top$ for a shortest preimage $u \in\Lambda$
of $\lambda_r \in \Lambda / 2 \Lambda$.
& \ref{ssect:mapping_axis} \\
\Skip{
MOG &     
Miracle Octad Generator; a tool for calculations in 
$\mathcal{C}  \subset \setF_2^{24}$.   
& \ref{xxx}   \\
}                 
$N_0$ &
A maximal subgroup of $\setM$ of structure
$2^{2+11+22}.(M_{24} \times S_3)$.
&  \ref{ssect:N0}  \\     
$N_{x0}$ &       
Subgroup of structure   $2_+^{1+24}.2^{11}.M_{24}$ of $\setM$, 
with $G_{x0} \cap N_0 = N_{x0}$.
& \ref{ssect:N0}  \\                
${\Omega} $   &
The  element $(\tilde{\Omega},0)$ of the Parker loop $\mathcal{P}$.
& \ref{ssect:Parker}    \\    
$\tilde{\Omega} $
&  The set $\{0, \ldots, 23\}$ used for labelling the basis vectors of
$\Field_2^{24}$ and &  \\
& $\setR^{24}$. The power set $2^{\tilde{\Omega}} $ is identified
with $\Field_2^{24}$; 
and we have  $\mathcal{C} \subset \Field_2^{24}$.
& \ref{ssect:Golay}    \\    
$\mathcal{P}$ &
The Parker loop, any $d \in \mathcal{P}$ has the form
$(\tilde{d}, \mu)$, $\tilde{d} \in \mathcal{C}, \mu \in \Field_2$.
& \ref{ssect:Parker}  \\
$\pi, \pi', \pi''$ &
Standard automorphisms of the Parker Loop  $\mathcal{P}$ 
in $\AutStP$.
& \ref{ssect:Parker}  \\   
$Q_{x0}$ &
A normal subgroup  of structure   $2_+^{1+24}$ 
of the group $G_{x0}$. 
&   \ref{ssect:Leech:Qx0} \\ 
$\rho$ &
196884-dimensional representation of $\setM$ with matrix entries
in  $\setZ[\frac{1}{2}]$. \\
& We have  
$\rho = 1_x \oplus 196883_x  = 
300_x  \oplus 98280_x  \oplus (24_x \otimes 4096_x)$.
& \ref{sect:rep:M} \\
$\rho_p$ &
Representation $\rho$ of $\setM$ with matrix entries
taken modulo $p$, $p$ odd.
& \ref{sect:rep:M} \\
$S_k(v)$ &
Set of short vectors in  $\Lambda/2 \Lambda$ depending
on $v \in \rho_{15}$ and integer $k$.
& \ref{ssect:mapping_axis} \\
$\mbox{type}(\lambda_r)$ &
Type of a vector $\lambda_r$ in $\Lambda$. This is half the squared 
length of  $\lambda_r$. \\ 
&For $\lambda_r \in \Lambda / 2 \Lambda$ this is the type of the shortest
preimage of $\lambda_r$ in  $\Lambda$.
&  \ref{ssect:Leech:Qx0}  \\
$\theta$ &
Cocycle  of the Parker loop $\mathcal{P}$, with
$(\tilde{d},0) \cdot (\tilde{e},0)
= (\tilde{d}+\tilde{e},\theta(\tilde{d},\tilde{e}))$.
&  \ref{ssect:Parker}   \\ 
$\tau$ &
Triality element, a generator of $\setM$ in $\Gamma$, which is in
$N_0 \setminus N_{x0}$.
&   \ref{ssect:N0} \\ 
$U(v)$ &
Subset of $\Lambda / 2\Lambda$ depending on an axis $v$
(and on the group $G_{x0}$).
&  \ref{ssect:mapping_axis} \\  
$U_4(v)$ &
Subset of $U(v)$ containing the vectors in $U(v)$ of type 4.
&  \ref{ssect:mapping_axis} \\  
$v_1$ &
A vector in $\rho_{15}$ fixed by the neutral element of $\setM$ only.
Any \\
&  unknown $g \in G_{x0}$ can be constructed 
from $v_1 g$ as a word in $\Gamma^*$.
& \ref{sect:Recogniz:Gx0} \\
$v^+, v^-$ &
Axes of the 2A involutions $x_\beta, x_{-\beta}$ 
in the representation $\rho$.
& \ref{sect:axes:2A} \\
\end{tabular}

\begin{tabular}{cll}
Symbol & Description 
\hphantom{01234567890123456789012345678901234567890123456789}
& \hspace{-2.25em} Section \vspace{0.5ex} \\
$x_{-1}$ &
The unique central involution in group $Q_{x0}$, 
and in group $G_{x0}$.
&  \ref{ssect:Leech:Qx0} \\
$x_\beta, x_{-\beta}$ &
Fixed involutions in $Q_{x0}$, both in class 2A in $\setM$;
$\lambda(x_{\pm\beta}) = \lambda_{\beta}$.
& \ref{ssect:central:H} \\  
$x_d, x_\delta, x_\pi$ &
Generators of $N_{x0}$	in $\Gamma$, for
$d \in \mathcal{P}, \delta \in \mathcal{C^*}, \pi \in \AutStP$.
&  \ref{ssect:N0}  \\      
$x_\Omega, x_{-\Omega}$ &
Fixed involutions in the group $Q_{x0}$, with
$\lambda(x_{\pm\Omega}) = \lambda_{\Omega}$. 
& \ref{ssect:N0} \\
$\xi$ &
A generator of $\setM$ in $\Gamma$, which is in
$G_{x0} \setminus N_{x0}$.
&   \ref{sect:rep:M} \\ 
$y_d$&
A generator of $N_{x0}$	in $\Gamma$, for $d \in \mathcal{P}$.
&  \ref{ssect:N0}  \\      
$1_\rho$ &
The unit matrix in the representation $300_x$; we have  
$1_\rho \in 1_x$.
& \ref{sect:rep:M}     \\
$1_x$ &
Trivial representation of $\setM$ or of $G_{x0}$;
we have $1_x \subset 300_x$.
&  \ref{sect:rep:M} \\     
$24_x$ &
Natural representation of $\mbox{Co}_0$ as the automorphism
group of $\Lambda$. &
  \ref{sect:rep:M} \\ 
$300_x$ &
Rational representation of $G_{x0}$, subspace of  $\rho$,
isomorphic to the \\
&  space of real symmetric $24 \times 24$ matrices; 
we have $\Lambda \subset \setR^{24}. $ 
& \ref{sect:rep:M} \\  
$4096_x$ &
A representation such that
$24_x \otimes 4096_x$ is a representation of $G_{x0}$.
&  \ref{ssect:maximal:G0} \\
$98280_x$ &
A monomial rational representation of $G_{x0}$, subspace of $\rho$.
&  \ref{sect:rep:M} \\
$196883_x$ &
Minimal real faithful representation of $\setM$, 
with $\rho =  1_x \oplus 196883_x$.
& \ref{sect:rep:M} \\
$\tilde{d}$ &
Image of $d \in \mathcal{P}$ in $\mathcal{C}$ under
the natural homomorphism $\mathcal{P} \rightarrow \mathcal{C}$.
& \ref{ssect:Parker}  \\  
$|S|, |\delta|$ &
Cardinality of a finite set $S$; minimum weight of
cocode word $\delta$.
& \ref{ssect:Golay}    \\
$\|M\|$ &
Sum of the squares of the entries of a symmetric matrix $M$.
& \ref{ssect:Enumeration:Gx0} \\ 
$\left<.,.\right>$ & 
The scalar product, e.g. on $\mathcal{C} \times \mathcal{C^*}$,
on $\Lambda \times \Lambda$, or on $\rho \times \rho$.
&  \ref{ssect:Golay}     \\
\end{tabular}

%%%%%%%%%%%%%%%%%%%%%%%%%%%%%%%%%%%%%%%%%%%%%%%%%%%%%%%%%%%%%%%%%%%%%%%%%%%%%%
%%%%%%%%%%%%%%%%%%%%%%%%%%%%%%%%%%%%%%%%%%%%%%%%%%%%%%%%%%%%%%%%%%%%%%%%%%%%%%
%%%%%%%%%%%%%%%%%%%%%%%%%%%%%%%%%%%%%%%%%%%%%%%%%%%%%%%%%%%%%%%%%%%%%%%%%%%%%%
\appendix
%%%%%%%%%%%%%%%%%%%%%%%%%%%%%%%%%%%%%%%%%%%%%%%%%%%%%%%%%%%%%%%%%%%%%%%%%%%%%%
%%%%%%%%%%%%%%%%%%%%%%%%%%%%%%%%%%%%%%%%%%%%%%%%%%%%%%%%%%%%%%%%%%%%%%%%%%%%%%
%%%%%%%%%%%%%%%%%%%%%%%%%%%%%%%%%%%%%%%%%%%%%%%%%%%%%%%%%%%%%%%%%%%%%%%%%%%%%%

%%%%%%%%%%%%%%%%%%%%%%%%%%%%%%%%%%%%%%%%%%%%%%%%%%%%%%%%%%%%%%%%%%%%%%%%%%%%%%
\section{Mapping a type-4 vector in  $\Lambda/ 2\Lambda$ to the standard frame}
\label{app:map:type4}
%%%%%%%%%%%%%%%%%%%%%%%%%%%%%%%%%%%%%%%%%%%%%%%%%%%%%%%%%%%%%%%%%%%%%%%%%%%%%%

In this appendix we show how to find an element $g$ of the group 
$G_{x0}$ defined in Section~\ref{ssect:maximal:G0} that maps 
an arbitrary type-4 vector $\lambda_r$ in $\Lambda / 2\Lambda$ 
to the standard frame $\lambda_\Omega$ in  $\Lambda / 2\Lambda$. 
The  normal subgroup $Q_{x0}$ of $G_{x0}$ of structure $2_+^{1+24}$ 
operates trivially on $\Lambda / 2\Lambda$.
Hence the group  $G_{x0}$ of structure $2_+^{1+24}.\mbox{Co}_1$ 
acts on  $\Lambda / 2\Lambda$ in the same way as $\mbox{Co}_1$ 
acts as automorphism group on $\Lambda / 2\Lambda$. So we may work 
in the factor group $\mbox{Co}_1$ of $G_{x0}$; and we may assume 
that $\mbox{Co}_1$ is generated by $y_d, x_\pi, \xi$, as defined in
Section~\ref{sect:N0}--\ref{sect:rep:M}. Taking these generators
modulo $Q_{x0}$, we may assume
$d \in \mathcal{C} / \{1, \tilde{\Omega}\}, \pi \in M_{24}$.
So given a $\lambda_r \in \Lambda / 2\Lambda$ of type 4, it suffices 
to find an element of $\mbox{Co}_1$ (represented as a word in the 
generators $y_d, x_\pi,$ and $\xi$) that maps 
$\lambda_r$ to $\lambda_\Omega$.

Modulo $Q_{x0}$, the group generated by $y_d, x_\pi$ is a maximal
subgroup $\bar{N}_{x0} = N_{x0} / Q_{x0}$ of 
$\mbox{Co}_1$ of structure $2^{11}.M_{24}$. From the discussion in
Section~\ref{ssect:maximal:G0} we see that $\bar{N}_{x0}$ is the 
stabilizer of $\lambda_\Omega$ in $\mbox{Co}_1$.
The action of $\mbox{Co}_1$ (or of $\bar{N}_{x0}$) on the Leech 
lattice $\Lambda$ is defined up to sign. 

In the sequel we will describe the orbits of $\bar{N}_{x0}$ on
$\Lambda / 2 \Lambda$. 

The orbits of the group $\bar{N}_{x0}$ on the vectors of
type 2, 3, and 4 on $\Lambda$ have been described in 
\cite{ivanov_1999},  Lemma 4.4.1. $\bar{N}_{x0}$ acts monomially on
the Leech lattice $\Lambda$. Thus an orbit $\lambda_r \bar{N}_{x0}$
of $\bar{N}_{x0}$ on $\Lambda / 2 \Lambda$ can be described by the 
{\em shapes} of the vectors in $\Lambda(\lambda_r)$ for any  
$\lambda_r \in \lambda_r \bar{N}_{x0}$,  where $\Lambda(\lambda_r)$
is the set of the shortest preimages of $\lambda_r$ in the Leech
lattice as in Section~\ref{ssect:Leech:Qx0}. Here the shape of a 
vector is the multiset of the absolute values of the co-ordinates of 
the vector. E.g. a vector of shape $(3^5 1^{19})$ has 5 co-ordinates 
with absolute value~3 and 19 co-ordinates with absolute value~1.

A vector of type 2 or 3 in $\Lambda / 2 \Lambda$ has two
opposite preimages of the same type in $\Lambda$; so 
its shape is uniquely defined.  A vector of type 4 in 
$\Lambda / 2 \Lambda$ has $2 \cdot 24$ preimages 
of type 4 in $\Lambda$ which are orthogonal except when equal 
or opposite; see e.g. \cite{Conway-SPLG}, \cite{ivanov_1999}. 
\Skip{
It is well known that a type-4 vector in $\Lambda / 2 \Lambda$ 
corresponds to a co-ordinate frame in the Leech lattice in 
standard co-ordinates. 
}

The table at Lemma 4.4.1 in \cite{ivanov_1999} assigns a name 
and a shape to each orbit of  $\bar{N}_{x0}$ on the vectors of type
2, 3, and 4 in  $\Lambda$. The table at  Lemma 4.6.1
in \cite{ivanov_1999} assigns a name and one or more shapes 
to each orbit of $\bar{N}_{x0}$ on the vectors of type 4 in  
$\Lambda / 2 \Lambda$. We reproduce this information for the 
orbits of $\bar{N}_{x0}$ on $\Lambda / 2 \Lambda$ in 
Table~\ref{table:orbits:Nx0}. In column  {\em Subtype} of that 
table we also assign a subtype  (which is a 2-digit number) 
to each orbit. 
The first digit of the subtype specifies the type of the vectors 
in that orbit and the second digit is used to distinguish between 
orbits of the same type. Hints for memorizing the  second digit 
are given in Table~\ref{table:subtypes:Nx0}.

For $d \in \mathcal{P}, \delta \in \mathcal{C}^*$ let 
$x_d, x_\delta \in Q_{x0}$ as in Section~\ref{ssect:N0}. 
Put $\lambda_d = \lambda(x_d)\in  \Lambda / 2 \Lambda$,
$\lambda_\delta = \lambda{(x_\delta)} \in  \Lambda / 2 \Lambda$,
as in Section \ref{ssect:Leech:Qx0}.
Then $\lambda_d$ is well defined also for $d \in \mathcal{C}$,
and each vector  $\lambda_r \in \Lambda / 2 \Lambda$ has a unique 
decomposition 
$\lambda_r = \lambda_d + \lambda_\delta, d \in \mathcal{C}, 
\delta \in \mathcal{C}^*$. The subtype of a vector 
$\lambda_d + \lambda_\delta \in \Lambda / 2\Lambda$ can
be computed from $d$ and $\delta$,
as indicated in the Table~\ref{table:orbits:Nx0}.

\begin{table} % [H]
\[
\begin{array}{|c|c|c|c|c|c|c|}
	\hline 
	\mbox{Sub-} & \mbox{Name} & \mbox{Shape}  & |d| & |\delta| &
	\langle d , \delta\rangle & \mbox{Remark} \\
	\mbox{type}  & \mbox{in \cite{ivanov_1999}} & & & & & \\
	\hline
	00 & 
	& (0^{24})  & 
	0 & 0  & 0  & \\     
	\hline
	20 & 
	\Lambda_2^4 & (4^2 0^{22}) & 
	0, 24  &  2  & 0 &    \\     
	\hline
	21 & 
	\Lambda_2^3 & (3^1  1^{23}) & 
	\mbox{any} & 1 & |d| / 4 &  \\     
	\hline
	22 & 
	\Lambda_2^2 & (2^8 0^{16}) & 
	8, 16  &  \mbox{even} & 0 &  1. \\     
	\hline
	31 & 
	\Lambda_3^5 & (5^1 1^{23}) & 
	\mbox{any} & 1 & |d| / 4 + 1 &  \\     
	\hline
	33 & 
	\Lambda_3^3 & (3^3 1^{21}) & 
	\mbox{any} & 3 & |d| / 4 &  \\     
	\hline
	34 & 
	\Lambda_3^4 & (4^1 2^{8} 0^{15}) & 
	8, 16 & \mbox{even}  &  1 &  \\    
	\hline
	36 & 
	\Lambda_3^2 & (2^{12} 0^{12}) & 
	12 & \mbox{even}  &  0 &  \\    
	\hline
	40 & 
	\bar{\Lambda}_4^{4a} & (4^4 0^{20}) & 
	0, 24 & 4 &  0 &  \\    
	\hline 
	42 & 
	\bar{\Lambda}_4^{6} & (6^1 2^7 0^{16}), (2^{16} 0^8) & 
	8, 16 & \mbox{even} &   0 & 2.\\    
	\hline
	43 & 
	\bar{\Lambda}_4^{5} & (5^1 3^2 1^{21}), (3^{5} 1^{19}) & 
	\mbox{any} & 3 & |d| / 4 + 1 &  \\    
	\hline
	44 & 
	\bar{\Lambda}_4^{4b} & (4^2 2^8 0^{14}), (2^{16} 0^8) & 
	8, 16 & \mbox{even} &   0 & 3.\\    
	\hline 
	46 & 
	\bar{\Lambda}_4^{4c} & (4^1 2^{12} 0^{11}) & 
	12 & \mbox{even} &   1 &  \\    
	\hline 
	48 & 
	\bar{\Lambda}_4^{8} & (8^1  0^{23}) & 
	24 & 0 &  0 &  \\    
	\hline 
\end{array} 
\] \\
\noindent\makebox[\textwidth][c]
{
\begin{minipage}[t]{0.7\textwidth}
Remarks: 
\vspace{-1ex}
\begin{enumerate}
	\item	
	$|\delta|/2 = 1 + |d|/8 \pmod{2}$, 
	$\delta' \subset d \Omega^{1 + |d|/8}$ for a suitable \\
	representative $\delta'$ of the cocode element $\delta$
	in $\setF_2^{24}$.
	\item
    \vspace{-1ex}
	$|\delta|/2 = |d|/8 \pmod{2}$, 
	$\delta' \subset d \Omega^{1 + |d|/8}$ for a suitable \\
	representative $\delta'$  of the cocode element $\delta$
	in $\setF_2^{24}$.
	\item
    \vspace{-1ex}
	None of the conditions stated in Remarks 1 and 2 holds.
\end{enumerate}
\end{minipage}
}
\caption{Orbits of $N_{x0}$ of $\Lambda/2\Lambda$}
\label{table:orbits:Nx0}    
\end{table}

Columns  {\em Name} and {\em Shape} in Table \ref{table:orbits:Nx0}
list the names and the shapes of the orbits as given in 
\cite{ivanov_1999}, Lemma 4.1.1 and 4.6.1. Columns $|d|$ and 
$|\delta|$ list conditions on the weight of a Golay code word $d$ and 
of (a shortest representative of) the Golay cocode element $\delta$, 
respectively. Column $\langle d, \delta \rangle$ lists 
conditions on the scalar product of  $d$ and  
$\delta$. 
%All this information can easily be derived from  
%\cite{ivanov_1999} and  \cite{Conway:Construct:Monster} 
%(or from \cite{Seysen20}).  
The information in the four rightmost columns of 
Table~\ref{table:orbits:Nx0} can be derived by using
(\ref{eqn:lambda:Leech}).
Conditions for vectors of type~2 are also given in
\cite{Conway:Construct:Monster}.

\begin{table} % [H]
\makebox[\textwidth][c]{	
\begin{minipage}[t]{0.7\textwidth}
\begin{tabular}{|c|l|}
	 \hline
	 Subtype & Meaning\\
	  \hline
	  x0 & Orbit contains a vector $\lambda_d+\lambda_\delta$ 
	      with $|d|=0$; $|\delta|$ even. \\
	 \hline
	  x1 & Vectors $\lambda_d+\lambda_\delta$ with $|\delta|=1$. \\
	 \hline
	  x2 & Vectors $\lambda_d+\lambda_\delta$ with $|d| = 8, 16$;
	      $|\delta|$ even;\\
	  &  and Remark 1 or 2 in Table~\ref{table:orbits:Nx0} 
	  applies. \\
	 \hline
	  x3 & Vectors $\lambda_d+\lambda_\delta$ with $|\delta|=3$. \\
	 \hline
	  x4 & Vectors $\lambda_d+\lambda_\delta$ with $|d| = 8, 16$;
	      $|\delta|$ even; \\
	     &  and none of the Remarks 1 or 2 in
	      Table~\ref{table:orbits:Nx0} apply. \\
	 \hline
	  x6 & Vectors $\lambda_d+\lambda_\delta$ with $|d| = 12$; 
	      $|\delta|$ even. \\
	 \hline
	  x8 & Reserved for the standard frame $\lambda_\Omega$. \\
	 \hline
\end{tabular}
\end{minipage}
}
\caption{Subtypes of orbits of $N_{x0}$ of $\Lambda/2\Lambda$}
\label{table:subtypes:Nx0}    
\end{table}

%Here the quantifiers run over vectors in $\mathbb{F}_2^{24}$ equal to
%$\delta \pmod{\mathcal{C}}$.

Table \ref{table:orbits:Nx0}   provides the information required
for effectively computing
the subtype of an element $\lambda_d + \lambda_\delta$ from
$d$ and $\delta$. 

From the generators $y_d, x_\pi,$ and $\xi$ of $\mbox{Co}_1$ 
mentioned above, only the generators $\xi^{\pm 1}$ may
change the subtype of a vector in $\Lambda/ 2 \Lambda$.
We say that a vector $v \in \mathbb{Z}^n$ is of shape
$(m^\alpha 0^{n-\alpha} \bmod 2m)$ if $v$ has
$\alpha$ co-ordinates equal to $m$ $\pmod{2m}$,
and $n-\alpha$ co-ordinates equal to 
$0$ $\pmod{2m}$.
For a vector $v = (v_0,\ldots, v_{23}) \in \mathbb{R}^{24}$
define 
$(v_{4i}, v_{4i+1}, v_{4i+2}, v_{4i+3}) \in \mathbb{R}^4$ 
to be the $i$-th column of $v$. 
This definition is related to the Miracle Octad Generator (MOG) 
used for the description of the Golay code, see 
\cite{Conway-SPLG}, Ch. 11. The co-ordinates of a vector in the
Leech lattice $\Lambda$ may be labelled with the entries of the 
MOG. Using this MOG labelling, we will depict the co-ordinates of
a vector $v \in \Lambda \subset \mathbb{Z}_2^{24}$  as follows:
\[
v = (v_0, \ldots, v_{23}) = 
\begin{array}{|c|c|c|c|c|c|}
	\hline 
     v_{0} &  v_{4} &  v_{8} & v_{12} & v_{16} & v_{20} \\
	\hline 
     v_{1} &  v_{5} &  v_{9} & v_{13} & v_{17} & v_{21} \\
	\hline 
     v_{2} &  v_{6} &  v_{10} & v_{14} & v_{18} & v_{22} \\
	\hline 
     v_{3} &  v_{7} &  v_{11} & v_{15} & v_{19} & v_{23} \\
	\hline 
\end{array} \quad .
\]
The following lemma provides some more information about the 
operation of  $\xi^k$ on $\Lambda$.

\begin{Lemma}
\label{Lemma:xi:MOG}
Let $v \in \Lambda \subset \mathbb{Z}^{24}$, and let 
$w$ be a column of $v$ in the MOG. Let $w^{(k)}$
be the corresponding column of $ v \cdot {\xi^k}$. Then
$w^{(k)}$ depends on $w$ only, and the squared
sums of the entries of $w$ and $w^{(k)}$ are equal.
If $w$ has shape
$(m^4 \bmod{2m})$ then there is a unique $k \in \{\pm1\}$
such that $w^{(k)}$ has shape $(0^4 \bmod{2m})$.
If $w$ has shape $(m^2 0^2 \bmod{2m})$ then 
$w^{(k)}$ has shape $(m^2 0^2 \bmod{2m})$ for
$k = \pm 1$.
\end{Lemma}

\fatline{Sketch proof}

In Equation (\ref{eqn:matrix:xi}) the operation of $\xi$ on a column 
of the MOG is given as a product of two orthogonal symmetric
$4 \times 4$ matrices $A$ and $B$. Thus $A$ and $B$ are involutions,
and $\xi^{-1}$ acts as multiplication by $BA$ on a column. Note that
multiplication of a column with $-A$ is equivalent to
subtracting the halved sum of the column from each entry, and
multiplication with $B$ means negating the first entry.
So it is easy to show the stated properties of the action
of $\xi^{\pm 1}$ on columns of shape  $(m^4 \bmod{2m})$
and $(m^2 0^2 \bmod{2m})$.

\proofend

We apply Lemma~\ref{Lemma:xi:MOG} to vectors $v \in \Lambda$ of 
type 4. Here the entries of the MOG columns of $v$ are either all
even or all odd; and they are usually quite small. So the shape
of a column of  $v$ imposes rather severe restrictions on the
shapes of the corresponding column of $v \xi^{\pm}$.
In the cases discussed below this allows to infer the 
subtype of  $v \xi$ or  $v \xi^{-1}$ from the shapes of the 
colums of $v$.

\begin{figure} [!h]
	 \centering
	\begin{tikzpicture}
    nodes={draw, % General options for all nodes
	line width=1pt,
	anchor=center, 
	text centered,
	rounded corners,
	minimum width=1.5cm, minimum height=8mm
    }, 
		\matrix (m) [matrix of  nodes,row sep=2.0ex,column sep=2em,
		  minimum width=2em, nodes={rounded corners=6pt, draw}]
		{   
			43 &    &  42   &     \\
			&      &       &   40 &  48  \\
		    46 	&      &  44 \\
		};         	
		\path[->] 
		(m-1-1) edge (m-1-3) 
	%	(m-1-1) edge (m-3-3) 
		(m-3-1) edge (m-3-3) 
		(m-1-3) edge (m-2-4) 
		(m-3-3) edge (m-2-4) 
		(m-2-4) edge (m-2-5) 
		;
		) ;
	\end{tikzpicture} 
	\caption{Reduction of the orbits of $\bar{N}_{x0}$ on type-4 vectors  
	  in $\Lambda/2\Lambda$ }  
	\label{figure:M24:orbits_type4}    
\end{figure}
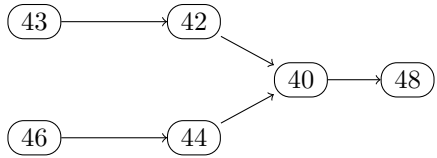

In the sequel we apply an operation $x_\pi \xi^k$ to a vector
$\lambda_r$ in $\Lambda / 2\Lambda$ of a given subtype 4X, so that 
the subtype is changed as in  Figure~\ref{figure:M24:orbits_type4}. 
The graph in that figure is a subgraph of the {\em Leech graph}
in \cite{ivanov_1999}, Section 4.7.
By a sequence of such operations we may eventually map an 
arbitrary vector in $\Lambda/2\Lambda$ to the (unique) vector 
$\lambda_\Omega$ of subtype 48. For the following discussion
of the individual subtypes we assume that the reader is familiar 
with the Mathieu group $M_{24}$, the Golay Code $\mathcal{C}$,
and the MOG, as  presented in  \cite{Conway-SPLG}, Ch. 11.

In the following discussion the phrase {\em up to MOG permutations}
means up to permutations of the columns of the MOG and  
arbitrary permutations inside a column of the MOG.

\vspace{1ex}
\fatline{From subtype 46 to subtype 44}
\nopagebreak

Let $\lambda_r = \lambda_d + \lambda_\delta$ be of 
subtype 46. Then $d$ is a dodecad. We first assume that $d$
contains one column of the MOG. 

From the discussion of the MOG in  \cite{Conway-SPLG} we
see that the union of two different columns of the MOG is
an octad, and that the intersections of $d$ with all columns 
must have the same size (modulo 2). Since a dodecad may not 
contain an octad, $d$ must intersect four columns of the MOG 
in a set of size 2; and its intersection with the 
remaining column is empty.

Note that each vector in $\Lambda(\lambda_r)$ has exactly one 
co-ordinate $\pm4$; and such a co-ordinate may occur at all 
24~positions, when considering all 48 vectors in
$\Lambda(\lambda_r)$.

Thus up to signs and MOG permutations there is a vector 
$v \in \Lambda(\lambda_r)$ with co-ordinates
\[
\begin{array}{|c|c|c|c|c|c|}
	\hline 
	2 & 0 & 2 & 0 & 0 & 0 \\
	\hline 
	2 & 0 & 2 & 0 & 0 & 0  \\ 
	\hline 
	2 & 0 & 4 & 2 & 2 & 2 \\ 
	\hline 
	2 & 0 & 0 & 2 & 2 & 2  \\ 
	\hline 
\end{array} \quad .
\]

By Lemma~\ref{Lemma:xi:MOG}  there is a 
$k = \pm 1$ such that column 0 of the vector 
$w = v \cdot {\xi^k}$ has shape $(4 \, 0^3)$. The other
columns of $w$ have the same shape as the corresponding 
columns of $v$. Thus  $w$ has shape 
$(4^1 2^8 0^{15})$. So from Table~\ref{table:orbits:Nx0}    
we see that  $w + 2\Lambda$ has subtype 44. 

\Skip{
Note that both, the dodecad $d$ and its complement, 
contain exactly one column of the MOG. 

The union of these two
columns is a grey even octad $o$ (in the notation of
\cite{Seysen20}).}

If dodecad  $d$ does not contain a column of the MOG then we 
select a permutation in $M_{24}$ that maps four entries of 
dodecad $d$ to the subset  $\{0,1,2,3\}$ of $\tilde{\Omega}$,
i.e. to the first column of the MOG. This is possible, since 
$M_{24}$ is quintuply transitive. 
Then we proceed as above.

\vspace{1ex}
\fatline{From subtype 43 to subtype 42}
\nopagebreak

Let $\lambda_r = \lambda_d + \lambda_\delta$ be of subtype 43.
Then $|\delta| = 3$. The three elements of (the shortest representative 
of) $\delta$ can lie in one, two, or three different columns of the MOG.

If all entries of  $\delta$ lie in the same column of the MOG
then up to signs and MOG permutations, there is a vector 
$v \in \Lambda(\lambda_r)$ with  co-ordinates
\[
\begin{array}{|c|c|c|c|c|c|}
	\hline 
	5 & 1 & 1 & 1 & 1 & 1 \\
	\hline 
	3 & 1 & 1 & 1 & 1 & 1  \\ 
	\hline 
	3 & 1 & 1 & 1 & 1 & 1 \\ 
	\hline 
	1 & 1 & 1 & 1 & 1 & 1  \\ 
	\hline 
\end{array}  \quad .
\]

By Lemma ~\ref{Lemma:xi:MOG} there is a $k = \pm 1$ such that 
one column of $w = v \cdot {\xi^k}$ (and hence all columns) 
have even entries. Column 0 of $w$ has squared norm 
$44 = 6^2 + 2^2 + 2^2 + 0^2$. That decomposition of
$44$ into a sum of four even squares is unique, so 
column 0 has shape $(6^1 2^2 0^1)$. The other columns
have shape  $(2^1 0^3)$. So $w$ is of
shape  $(6^1 2^7 0^{15})$; and hence  
$\lambda_r \xi^k = w + 2 \Lambda$ is of subtype 42.

If the entries of  $\delta$ lie in two or three different columns 
of  the MOG then we apply a permutation in $M_{24}$ that maps
$\delta$ to  $\{0,1,2\}$, and proceed as above.

\vspace{1ex}
\fatline{From subtype 44 to subtype 40}
\nopagebreak

Let $\lambda_r = \lambda_d + \lambda_\delta$ be of subtype 44. 
Then $d$ is an octad or a complement of
an octad. We call that octad $o$.
If the cocode word $\delta$ is a duad (i.e. $|\delta| = 2$) then
we have  $\delta = \{c_0, c_1\}$ and $o \cap \{c_0, c_1\} = \{\}$.
Otherwise, $\delta$ 
corresponds to a sextet containing a tetrad $c$ intersecting with 
$o$ in two points; in that case we let $c_0, c_1$ be the two 
elements of $c$ not contained in $o$.

Assume first that $o$ is a  union of
two columns of the MOG, and that the points  $c_0, c_1$
are in the same column of the MOG.

Then up to signs and MOG permutations there is a vector 
$v \in \Lambda(\lambda_r)$ with co-ordinates
\[
\begin{array}{|c|c|c|c|c|c|}
	\hline 
	2 & 2 & 4 & 0 & 0 & 0 \\
	\hline 
	2 & 2 & 4 & 0 & 0 & 0  \\ 
	\hline 
	2 & 2 & 0 & 0 & 0 & 0 \\ 
	\hline 
	2 & 2 & 0 & 0 & 0 & 0  \\ 
	\hline 
\end{array} \quad .
\]

Similar to the previous cases, we can show that there is a 
$k = \pm 1$ such that $w = v \cdot {\xi^k}$ has 
shape $(4^4 0^{20})$. So $w$ is of subtype 40.

Assume now that $o$ is any octad and that none of the two points  
$c_0, c_1 \in \tilde{\Omega}$ is in $o$. In the sequel we will
construct a permutation $\pi \in M_{24}$ that maps $o$ to a 
union of two columns of the MOG and $(c_0, c_1)$ to $(2,3)$. 
Then we may proceed with $\pi(o)$ as above. 

For constructing $\pi$ we select three elements $o_0, o_1, o_3$
from $o$; and we compute the syndrome $\sigma$  of the set 
$\{c_0, c_1, o_0, o_1, o_2\}$, which is a cocode word of length~3.
$\sigma$ intersects with  $o$ in exactly one point
$o_3$. Then we construct a permutation in $M_{24}$ that 
maps the tuple $(c_0, c_1, o_0, o_1, o_2, o_3)$ to the tuple
$(2, 3, 4, 5, 6, 7)$. Such a permutation exists, since both
tuples are subsets of an octad, $M_{24}$ is transitive on octads,
and the stabilizer of an octad in $M_{24}$ acts as the alternating 
permutation group on that octad. See \cite{Conway-SPLG} for proofs
of these statements.
This permutation maps $\{c_0, c_1\}$ to  $\{2,3\}$; and it maps
$o$ to an octad containing the column  $\{4,5,6,7\}$ of the MOG. 
Any octad containing one column of the MOG consists of two 
columns of the MOG.

\vspace{1ex}
\fatline{From subtype 42 to subtype 40}
\nopagebreak

Let $\lambda_r = \lambda_d + \lambda_\delta$ be of subtype 42. 
Then $d$ is an octad or a complement of an octad. We call that
octad $o$.

Assume first that $o$ is a  union of
two columns of the MOG.
Then up to signs and MOG permutations there is a vector 
$v \in \Lambda(\lambda_r)$ with co-ordinates
\[
\begin{array}{|c|c|c|c|c|c|}
	\hline 
	6 & 2 & 0 & 0 & 0 & 0 \\
	\hline 
	2 & 2 & 0 & 0 & 0 & 0  \\ 
	\hline 
	2 & 2 & 0 & 0 & 0 & 0 \\ 
	\hline 
	2 & 2 & 0 & 0 & 0 & 0  \\ 
	\hline 
\end{array} \quad .
\]

Similar to the previous cases, we can use 
Lemma~\ref{Lemma:xi:MOG} to show that there is a 
$k = \pm 1$ such that $w = v \cdot {\xi^k}$ has 
shape $(4^4 0^{20})$. So $w$ is of subtype 40.

Otherwise we first map the  octad  $o$ to the octad 
$\{0,1,2,3,4,5,6,7\}$, using a permutation in  $M_{24}$.

\vspace{1ex}
\fatline{From subtype 40 to subtype 48}
\nopagebreak

Let $\lambda_r$ be of subtype 40. Then
$\lambda_r = \alpha \lambda_\Omega + \lambda_\delta$,
$\alpha = 0, 1$, for some 
$\delta \in \mathcal{C}^*$, $|\delta|=4$.

If $\delta$ is equal to the standard tetrad  
represented by $\{0,1,2,3\}$, then a column of a
vector $v \in \Lambda(\lambda_r)$ in the MOG has shape $(4^4)$. 
By Lemma~\ref{Lemma:xi:MOG}  there is a 
$k = \pm 1$ such that $\xi^k$ maps this
column to a column of shape  $(8^1 0^3)$. 
Thus $\lambda_r \xi^k$ is of subtype 48 and hence equal to
$\lambda_{\Omega}$.

Otherwise we first apply a permutation in  $M_{24}$ that 
maps a tetrad in  $\delta$ to  $\{0,1,2,3\}$.

%%%%%%%%%%%%%%%%%%%%%%%%%%%%%%%%%%%%%%%%%%%%%%%%%%%%%%%%%%%%%%%%%%%%%%%%%%%%%%
\section{The probability that the kernel of a random symmetric matrix  
	on $\Lambda/ 3\Lambda$ is spanned by a type-4 vector}
\label{app:probability}
%%%%%%%%%%%%%%%%%%%%%%%%%%%%%%%%%%%%%%%%%%%%%%%%%%%%%%%%%%%%%%%%%%%%%%%%%%%%%%

We first calculate the probability $p^{(q)}_{n, k}$ that a symmetric 
$n \times n$ matrix with random entries in
$\mathbb{F}_q$ has corank $k$.  For any symmetric 
$n \times n$  matrix $Q_{n}$ of corank $k$  we can find
a basis of the underlying vector space such that $Q_n$ has shape
\[
Q_n =
\left(
\begin{array}{ccc}
	0    &  0      \\
	0    &  A   \\
\end{array}
\right) \, , \quad \mbox{where } A  \mbox{ is a non-singular symmetric } 
(n-k) \times (n-k) \mbox{ matrix. } 
\]
Let $Q_{n+1}$ be the symmetric $(n+1) \times (n+1)$ matrix obtained from
$Q_n$ by adding one row and one column:
\[
Q_{n+1} =
\left(
\begin{array}{ccc}
	0    &  0   &  b_0^\top \\
	0    &  A   &  b_1^\top   \\
	b_0  & b_1  &  c \\
\end{array}
\right)  \, , \quad \mbox{with } \quad
b_0  \in  \mathbb{F}_q^k, \;  b_1  \in  \mathbb{F}_q^{n-k}, \;
c \in \mathbb{F}_q \; .
\]
Matrix $Q_{n+1}$ has corank $k-1$ if $b_0 \neq 0$. In case $b_0 = 0$  
it has corank $k+1$ if $c = b_1 A^{-1} b_1^\top$ and corank $k$ otherwise.
Assuming that  $b_0, b_1,$ and $c$ are selected at random this implies:
\[
\mbox{corank}(Q_{n+1}) = \left\{
\begin{array}{c}
	k - 1   \\
	k      \\
	k + 1  \\
\end{array}
\right\}
\mbox{ with probability }
\left\{   
\begin{array}{cl}
	1 - q^{-k}   & \quad (\mbox{for } \, k > 0) \\
	q^{-k} - q^{-(k+1)}      \\
	q^{-(k+1)}  & \quad    \\
\end{array}
\right. \; .
\]
$Q_{n+1}$ has been  obtained from an arbitrary  $n \times n$ matrix $Q_n$ 
over $\mathbb{F}_q$ with corank~$k$. Hence for fixed $q$ the distribution
of $p^{(q)}_{n, k}, 0 \leq k \leq n,$ can be modelled as a Markov chain with
discrete time parameter $n$ and start value $p^{(q)}_{0, 0} = 1$.
So we easily obtain $p^{(3)}_{24, 1} \approx 0.31950$.

There are 398034000 vectors of type 4 in the Leech lattice  $\Lambda$,
see e.g. \cite{Conway-SPLG} Ch.4.11, and 
$3^{24}-1$ nonzero vectors in  $\Lambda / 3 \Lambda$,  so a random 
one-dimensional subspace of  $\Lambda / 3 \Lambda$
is spanned by a type-4 vector with probability about $1/709.56$. Thus 
the kernel of a random symmetric $24 \times 24$  matrix over 
$\mathbb{F}_3$ (acting on $\Lambda / 3 \Lambda$) is spanned by a 
type-4 vector with probability about $1/2221$.

Let $M(.)$ be as in Section~\ref{sect:Recogniz:Gx0}.
For $v \in \rho_3$ let $M_3(v)$ be the matrix $M(v)$ with entries
taken modulo 3.
Generating a vector $v$ in $\rho_3$ that is fixed by an element of 
order 71 costs 70 group operations on  $\rho_3$. So in average we need 
about 155000 group operations on $\rho_3$ to find a such vector 
$v \in \rho_3$ with the additional property that $\ker M_3(v)$
is spanned by a type-4 vector (modulo  $3\Lambda$). This is certainly 
doable, and has been done. Once having found a vector $v$ fixing an
element of order 71,  we may also check $\ker M_3(v \pm 1_\rho)$.

So we can find a vector $v$ satisfying the properties mentioned above,
such that $\ker M_3(v)$ is spanned by $w + 3\Lambda$, for
a type-4 vector $w \in \Lambda$.  Given the co-ordinates of $w$ 
(modulo 3), we can easily compute the corresponding type-4 vector 
in $\Lambda$, and hence also the corresponding  type-4 vector 
$\lambda_r$ in $\Lambda / 2\Lambda$. We put $v_1 = v g$, where $g$ 
is an element of $G_{x0}$ that maps the type-4 vector $\lambda_r$ to
the standard frame $\lambda_\Omega$. We use the method in
Appendix~\ref{app:map:type4} for computing a suitable
$g \in G_{x0}$. Then $\ker M_3(v_1)$ is spanned by 
$\eta_j + 3\Lambda$ for some basis vector $\eta_j$ of $\setR^{24}$.

%%%%%%%%%%%%%%%%%%%%%%%%%%%%%%%%%%%%%%%%%%%%%%%%%%%%%%%%%%%%%%%%%%%%%%%%%%%%%%
\section{Mapping a type-2 vector in  $\Lambda/ 2\Lambda$ to the standard 
	type-2 vector $\lambda_\beta$}
\label{app:map:type2}
%%%%%%%%%%%%%%%%%%%%%%%%%%%%%%%%%%%%%%%%%%%%%%%%%%%%%%%%%%%%%%%%%%%%%%%%%%%%%%

In this appendix we show how to find an element $g$ of $G_{x0}$ that 
maps an arbitrary type-2 vector $\lambda_r$ in $\Lambda / 2\Lambda$ 
to the standard type-2 vector $\lambda_\beta$ in  $\Lambda / 2\Lambda$.

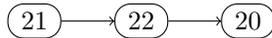
\begin{figure}  [!h]
	\centering
	\begin{tikzpicture}
		nodes={draw, % General options for all nodes
			line width=1pt,
			anchor=center, 
			text centered,
			rounded corners,
			minimum width=1.5cm, minimum height=8mm
		}, 
		\matrix (m) [matrix of  nodes,row sep=2.0ex,column sep=2em,
		minimum width=2em, nodes={rounded corners=6pt, draw}]
		{   
			21 &   22  &  20    \\
		};         	
		\path[->] 
		(m-1-2) edge (m-1-3) 
		(m-1-1) edge (m-1-2) 
		;
		) ;
	\end{tikzpicture} 
	\caption{Reduction of the orbits of $\bar{N}_{x0}$ on type-2 vectors  
		in $\Lambda/2\Lambda$ }  
	\label{figure:M24:orbits_type2}    
\end{figure}

We use the same notation as in Appendix~\ref{app:map:type4}.
We apply an operation $x_\pi \xi^k$ to a vector in
$\Lambda / 2\Lambda$ of a given subtype 2X, so that the subtype is 
changed as in  Figure~\ref{figure:M24:orbits_type2}. Finally,
we apply an operation in $N_{x0}$ to a vector of subtype 20
in order to map that vector to $\lambda_\beta$ .

\vspace{1ex}
\fatline{From subtype 21 to subtype 22}
\nopagebreak

Let $\lambda_r $ be of subtype 21.  Then $\Lambda(\lambda_r)$ is of 
shape $(3^1 1^{23})$. By Lemma~\ref{Lemma:xi:MOG} there is a element 
$\xi^k$ such that $\Lambda(\lambda_r) \xi^k$ has even co-ordinates. 
Then $\lambda_r \xi^k$ 
cannot be of subtype 21; and it is easy to see that  $\lambda_r \xi^k$ 
is of subtype~22.

\vspace{1ex}
\fatline{From subtype 22 to subtype 20}
\nopagebreak

Let $\lambda_r = \lambda_d + \lambda_\delta$ be of subtype 22. 
Then $d$ is an octad or a complement of an octad. We call that
octad $o$. Then  $\Lambda(\lambda_r)$ has shape $(2^8 0^{16})$.

Assume first that $o$ is a union of two columns of the MOG. 
Then by Lemma~\ref{Lemma:xi:MOG} there is a $k = \pm 1$ 
such that $\Lambda(\lambda_r) \cdot {\xi^k}$ has shape 
$(4^2 0^{22})$. So $\lambda_r \xi^k$ is of subtype 40.

Otherwise we first map the  octad  $o$ to the octad 
$\{0, 1, 2, 3, 4, 5, 6, 7\}$, using a permutation in  $M_{24}$.

\vspace{1ex}
\fatline{From subtype 20 to the vector  $\lambda_\beta$ }
\nopagebreak

Let $\lambda_r$ be of subtype 20.  Then
$\lambda_r = \alpha \lambda_\Omega + \lambda_\delta$,
$\alpha \in \setF_2$, for some 
$\delta \in \mathcal{C}^*$, $|\delta|=2$. So there is an 
$x_\pi$ that maps $\delta$ to $\beta$; i.e.
$\lambda_r x_\pi= \alpha' \lambda_\Omega + \lambda_\beta$.
In case $\alpha' = 0$ we are done. Otherwise we apply a
transformation $y_d$ with $\langle d, \beta \rangle = 1$.
So we have $\lambda_r x_\pi y_d^{\alpha'}=  \lambda_\beta$.

%%%%%%%%%%%%%%%%%%%%%%%%%%%%%%%%%%%%%%%%%%%%%%%%%%%%%%%%%%%%%%%%%%%%%%%%%%%%%%
\section{Mapping a feasible  type-2 vector in  $\Lambda/ 2\Lambda$ to the 
	type-2	vector $\lambda_{\Omega} \! + \! \lambda_\beta$, fixing the 
	standard type-2 vector $\lambda_\beta$}
\label{app:map:feasible}
%%%%%%%%%%%%%%%%%%%%%%%%%%%%%%%%%%%%%%%%%%%%%%%%%%%%%%%%%%%%%%%%%%%%%%%%%%%%%%

In this appendix we show how to find an element $g$ of $H$ that 
maps an arbitrary feasible type-2 vector 
$\lambda_r$ in $\Lambda / 2\Lambda$ to the standard feasible vector
$\lambda_\Omega + \lambda_\beta$ in  $\Lambda / 2\Lambda$.
Here the term {\em feasible} is defined as in 
Section~\ref{ssect:Enumeration:H}.

We use the same notation as in Appendix~\ref{app:map:type4}.
We apply an operation $x_\pi \xi^k$, with $\pi \in H \cap \AutStP$
to a feasible vector in
$\Lambda / 2\Lambda$ of a given subtype 2X, so that the subtype is 
changed as in  Figure~\ref{figure:M24:orbits_type2}. Note that
$H \cap \AutStP$ operates on set 
$\tilde{\Omega} = \{0, \ldots, 23\}$ as the
subgroup of the permutation group $M_{24}$ fixing the set  
$\beta = \{2,3\}$. Finally, we apply an operation in $H$  
to a feasible  vector of subtype 20 in order to map that vector
to $\lambda_\Omega + \lambda_\beta$.

\vspace{1ex}
\fatline{From subtype 21 to subtype 22}
\nopagebreak

The same argument as in the corresponding case in
Appendix~\ref{app:map:type2} shows that there is a
$k = \pm{1}$ such that $\lambda_r \xi^k$ is of subtype~22.
Since $\xi \in H$, the vector  $\lambda_r \xi^k$  is feasible.

\vspace{1ex}
\fatline{From subtype 22 to subtype 20}
\nopagebreak

Let $\lambda_r = \lambda_d + \lambda_\delta$ be of subtype 22. 
Then $d$ is an octad or a complement of an octad. We call that
octad $o$. Then  $\Lambda(\lambda_r)$ has shape $(2^8 0^{16})$.
Since $\lambda_r$ is feasible, the set $\{2,3\}$ is either 
contained in or disjoint to the  octad $o$; otherwise we 
would have $\mbox{type}(\lambda_r + \lambda_\beta) \neq 4$.

Assume first that $o$ is a  union of two columns of the MOG. 
Then by Lemma~\ref{Lemma:xi:MOG} there is a $k = \pm1$ 
such that $\Lambda(\lambda_r) \cdot {\xi^k}$ has shape 
$(4^2 0^{22})$. So $\lambda_r \xi^k$ is of subtype 20.

Otherwise, if the set $\{2,3\}$ is contained in $o$ then we first map 
the  octad  $o$ to the octad $\{0, 1, 2, 3, 4, 5, 6, 7\}$, using a
permutation in  $M_{24}$ that fixes the set $\{2,3\}$. Then
we proceed as above.

If the set $\{2,3\}$ is disjoint from $o$ then we construct 
a permutation $\pi \in M_{24}$ that maps octad $o$ to a union
of two columns of the MOG, and that fixes the tuple $(2,3)$.
Therefore we use the method in  Appendix~\ref{app:map:type4}
in the case 'From subtype 44 to subtype 40'. Then we proceed 
as above.

\vspace{1ex}
\fatline{From subtype 20 to the vector  $\Lambda_\Omega + \lambda_\beta$ }
\nopagebreak

Let $\lambda_r$ be feasible and of subtype 20.  Then
$\lambda_r = \alpha \lambda_\Omega + \lambda_\delta$,
$\alpha = 0, 1$, for some 
$\delta \in \mathcal{C}^*$, $|\delta|=2$. 
In case $\delta = \beta$  we have
$\lambda_r = \lambda_\Omega + \lambda_\beta$, so that we are done.

Otherwise both elements of $\delta$  are different from 2 and 3.
In case  $\delta = \{0,1\}$ we have 
$(-1)^{\alpha+1} \cdot 4 \eta_0 + 4 \eta_1 \in \Lambda(\lambda_r)$; 
and a direct calculation using (\ref{eqn:matrix:xi}) shows 
$\lambda_r \xi^{ 2 - \alpha} = \lambda_\Omega + \lambda_\beta$.

In case  $\delta \neq \{0,1\}$ we first apply a  permutation
in $M_{24}$ that maps $\delta$ to $\{0,1\}$ and fixes 
$\beta = \{2,3\}$. Then we proceed as in case  $\delta = \{0,1\}$.

%\vfill

\newpage

%\bibliography{../BibTex/references}{}

\begin{thebibliography}{10}
	
	\bibitem{Aschbacher-Sporadic}
	M.~Aschbacher.
	\newblock {\em Sporadic Groups}.
	\newblock Cambridge University Press, New York, 1986.
	
	\bibitem{Conway:Construct:Monster}
	J.~H. Conway.
	\newblock A simple construction of the {F}ischer-{G}riess monster group.
	\newblock {\em Inventiones Mathematicae}, 79, 1985.
	
	\bibitem{Atlas}
	J.~H.\ Conway, R.~T.\ Curtis, S.~P.\ Norton, R.~A.\ Parker, and R.~A.\ Wilson.
	\newblock {\em Atlas of Finite Groups}.
	\newblock Clarendon Press, Oxford, 1985.
	
	\bibitem{Conway-SPLG}
	J.~H. Conway and N.~J.~A. Sloane.
	\newblock {\em Sphere Packings, Lattices and Groups}.
	\newblock Springer-Verlag, New York, 3rd edition, 1999.
	
	\bibitem{DLP2023}
	H.~Dietrich, M.~Lee, and T.~Popiel.
	\newblock {The maximal subgroups of the Monster}.
	\newblock {\em arXiv e-prints}, page arXiv:2304.14646, May 2023.
	
	\bibitem{Far12}
	A.~Farooq.
	\newblock Some computations in the monster group and related topics.
	\newblock Doctoral thesis, School of Mathematical Sciences. Queen Mary
	University of London, 2012.
	
	\bibitem{Griess:Friendly:Giant}
	R.~L. Griess.
	\newblock The friendly giant.
	\newblock {\em Invent. Math}, 69:1--102, 1982.
	
	\bibitem{holmes_wilson_2003}
	P.~E. Holmes and R.~A. Wilson.
	\newblock A new computer construction of the monster using 2-local subgroups.
	\newblock {\em Journal of the London Mathematical Society}, 67(2):349--364,
	2003.
	
	\bibitem{ivanov_1999}
	A.~A. Ivanov.
	\newblock {\em Geometry of Sporadic Groups I, Peterson and Tilde Geometries}.
	\newblock Encyclopedia of Mathematics and its Applications. Cambridge
	University Press, 1999.
	
	\bibitem{citeulike:Monster:Majorana}
	A.~A. Ivanov.
	\newblock {\em The Monster Group and Majorana Involutions}.
	\newblock Cambridge University Press, April 2009.
	
	\bibitem{LPWW-Computer_Monster}
	S.~Linton, R.~Parker, P.~Walsh, and R.~Wilson.
	\newblock Computer construction of the monster.
	\newblock {\em Journal of Group Theory}, 1:307--337, 1998.
	
	\bibitem{mueller_2008}
	J.~M\"uller.
	\newblock On the action of the sporadic simple baby monster group on its
	conjugacy class $\mbox{2B}$.
	\newblock {\em LMS Journal of Computation and Mathematics}, 11:15–27, 2008.
	
	\bibitem{MNW_2006}
	J.~M\"uller, M.~Neunh\"offer, and R.A. Wilson.
	\newblock Enumerating big orbits and an application: $\mbox{B}$ acting on the
	cosets of $\mbox{Fi}_{23}$.
	\newblock {\em Journal of Algebra}, 314:75–96, 2006.
	
	\bibitem{Nor02}
	S.~Norton.
	\newblock Transforming the monster presentation.
	\newblock Preprint, published as an appendix in \cite{Far12}, 2002.
	
	\bibitem{Norton98}
	S.~P. Norton.
	\newblock Anatomy of the monster \uppercase{I}.
	\newblock In {\em The Atlas of Finite Groups: Ten Years On}, pages 198--214.
	Cambridge University Press, 1998.
	
	\bibitem{Seysen20}
	M.~Seysen.
	\newblock {A computer-friendly construction of the monster}.
	\newblock {\em arXiv e-prints}, page arXiv:2002.10921, February 2020.
	
	\bibitem{mmgroup2020}
	M.~Seysen.
	\newblock {\tt mmgroup}, a python implementation of the monster group.
	\newblock \url{https://github.com/Martin-Seysen/mmgroup}, 2020.
	\newblock Release {\tt mmgroup v1.0.0}, accessed on Jan 23, 2024.
	
	\bibitem{mmgroup_doc}
	M.~Seysen.
	\newblock Welcome to mmgroup's documentation.
	\newblock \url{https://mmgroup.readthedocs.io/en/stable}, 2021.
	\newblock Accessed on Jan 23, 2024.
	
	\bibitem{Wilson13}
	R.~A. Wilson.
	\newblock {The Monster and black-box groups}.
	\newblock {\em arXiv e-prints}, page arXiv:1310.5016, October 2013.
	
\end{thebibliography}
%\bibliographystyle{plain}

%\Skip{

%}

\end{document}